\documentclass[11pt]{amsart}
\usepackage[a4paper,margin=1.05in]{geometry}
\usepackage[T1]{fontenc}
\usepackage{lmodern,amsmath,amssymb,amsthm,mathtools,enumitem,microtype,needspace}
\usepackage[colorlinks,linkcolor=blue,citecolor=blue,urlcolor=blue,hypertexnames=false]{hyperref}
\newtheorem{theorem}{Theorem}[section]
\newtheorem{proposition}[theorem]{Proposition}
\newtheorem{lemma}[theorem]{Lemma}
\newtheorem{corollary}[theorem]{Corollary}
\theoremstyle{definition}

\newtheorem{remark}[theorem]{Remark}

\newcommand{\D}{\mathbb D}
\newcommand{\T}{\mathbb T}

\newcommand{\wh}{\widehat}
\newcommand{\one}{\mathbf 1}
\newcommand{\bh}{\mathbf h}
\newcommand{\bg}{\mathbf g}
\newcommand{\by}{\mathbf y}
\newcommand{\cE}{\mathcal E}
\newcommand{\cS}{\mathcal S}
\newcommand{\cQ}{\mathcal Q}
\newcommand{\E}{\mathbb E}
\newcommand{\BMOA}{\mathrm{BMOA}}

\title{Bounded N\"orlund operators on Hardy spaces}
\author{Georgios Stylogiannis}
\date{14 November 2026}
\begin{document}
\begin{abstract}
We characterize the nonnegative N\"orlund matrices that induce bounded
operators on Hardy spaces. For $a_n\ge0$ and $a_0>0$, boundedness on one
$H^p$, $1<p<\infty$, is equivalent to boundedness on every such space and
to a single quadratic testing condition involving the ratios of the
partial sums $A_n=\sum_{k=0}^na_k$. The testing constant is comparable
to the squared $H^2$ operator norm, with absolute constants. No doubling
or monotonicity assumption on the weights is required. The proof combines
an exact vector-valued Hankel representation with square-function
estimates and a geometric argument showing that the doubling-time
endpoints have finite lacunary order. Equivalently, boundedness is characterized by a
truncated quadratic block condition on the coefficients of
$G_a(z)=\sum_n(a_n/A_n)z^n$. Within the bounded family, $G_a$ belongs
to the Bloch space if and only if it belongs to $\BMOA$, and this occurs
exactly when $(A_n)$ is doubling.
\end{abstract}
\maketitle

\section{Introduction}

The Ces\`aro operator averages the Taylor coefficients of an analytic
function. A N\"orlund operator allows arbitrary nonnegative weights in
this average. Which weights give a bounded operator on $H^p$, and how
does the answer depend on $p$? We show that a single family of polynomial
tests answers both questions throughout the range $1<p<\infty$.

Let $a=(a_n)_{n\ge0}$ with $a_n\ge0$ and $a_0>0$, and put
$A_n=\sum_{j=0}^na_j$. The N\"orlund matrix acts on sequences
$s=(s_k)_{k\ge0}$ by
\begin{equation}\label{eq:matrix}
 (N_as)_n=\frac1{A_n}\sum_{k=0}^na_{n-k}s_k.
\end{equation}
Taking $s_k=\wh f(k)$ gives the analytic N\"orlund operator
\begin{equation}\label{eq:operator}
 \wh{N_af}(n)=\frac1{A_n}\sum_{k=0}^na_{n-k}\wh f(k).
\end{equation}
Each row is a convex average. For an analytic $f$ and any $0<r<1$,
Cauchy's estimate therefore gives
$|\wh{N_af}(n)|\le r^{-n}\max_{|z|=r}|f(z)|$, so the resulting series
defines an analytic function on $\D$. On $H^p=H^p(\D)$ we use the norm
\[
 \|f\|_{H^p}=\sup_{0<r<1}
 \left(\int_{\T}|f(r\zeta)|^p\,dm(\zeta)\right)^{1/p},
 \qquad 1<p<\infty,
\]
where $dm$ is normalized arc length. For $a_n=1$, \eqref{eq:matrix}
takes the arithmetic mean of $s_0,\ldots,s_n$, and \eqref{eq:operator}
gives the classical Ces\`aro operator
\[
 Cf(z)=\frac1z\int_0^z\frac{f(\zeta)}{1-\zeta}\,d\zeta,
 \qquad Cf(0)=f(0).
\]
It is bounded on $H^p$ for $1<p<\infty$, and its Hardy-space theory
is connected with composition semigroups \cite{Siskakis}.
More generally, $a_n=\binom{n+\alpha-1}{n}$, with $\alpha>0$, gives the
Ces\`aro method of order $\alpha$.

In classical summability theory, $s_k$ is the $k$th partial sum of a
series. The methods go back to Voronoi and N\"orlund \cite{Norlund}.
Hardy \cite[Chapter~IV]{Hardy} gives the classical treatment, and Hille
and Tamarkin \cite{HT} and Sahney \cite{Sahney} study the summability
of Fourier series. In this setting the first question is regularity,
or preservation of limits of convergent sequences. For nonnegative
weights, regularity is equivalent to $a_n/A_n\to0$
\cite[\S4.2, Theorem~16]{Hardy}. Boundedness is a
different requirement. The weights $a_n=e^{\sqrt n}$ give a regular
method whose matrix is unbounded on every $\ell^p$, $1<p<\infty$
\cite[Theorems~2 and~3]{BC}.

Borwein and Cass \cite{BC} obtained several sufficient conditions for
$\ell^p$ boundedness, including
\begin{equation}\label{eq:BC}
 \sup_{n\ge1}\frac{na_n}{A_n}<\infty.
\end{equation}
Cass and Kratz \cite{CK} obtained necessary and sufficient conditions
for weights generated by logarithmico-exponential functions. Borwein
\cite[Theorem~1(i)]{Bo93} proved that if $na_n/A_n$ has a limit in $[0,\infty]$,
then $\ell^p$ boundedness is equivalent to the finiteness of that limit.
He also gave bounded examples for which \eqref{eq:BC} fails.
By Parseval's identity, taking Taylor coefficients identifies $H^2$
isometrically with $\ell^2$. Under this identification, the action of
$N_a$ on $H^2$ corresponds to its matrix action on $\ell^2$ with the same
operator norm. For the other Hardy exponents, one must control the
boundary $L^p$ norm. The theorem below gives a criterion in terms of
$(A_n)$ alone and shows that it governs both problems.

\subsection{The characterization}

Test $N_a$ on $D_J(z)=1+z+\cdots+z^{J-1}$, where $J\ge1$.
Since $\|D_J\|_{H^2}^2=J$, the squared norm ratio is
\[
 \frac{\|N_aD_J\|_{H^2}^2}{\|D_J\|_{H^2}^2}
 =\frac{\|N_aD_J\|_{H^2}^2}{J}.
\]
With $A_j=0$ for $j<0$, the Taylor coefficients of $N_aD_J$ are
\begin{equation}\label{eq:rowsum}
 \wh{N_aD_J}(n)=R_J(n):=1-\frac{A_{n-J}}{A_n}.
\end{equation}
For $n\ge J$, $R_J(n)$ is the fraction of the mass of $A_n$ contributed
by the last $J$ weights. For $0\le n<J$, it equals one. Taking the
supremum of these ratios gives
\begin{equation}\label{eq:tau}
 \tau(a)=\sup_{J\ge1}\frac{\|N_aD_J\|_{H^2}^2}{J}
 =\sup_{J\ge1}\left[1+\frac1J\sum_{n\ge J}
       \left(1-\frac{A_{n-J}}{A_n}\right)^2\right],
\end{equation}
where either expression is allowed to be infinite.

\begin{theorem}\label{thm:main}
Let $a_n\ge0$ and $a_0>0$. The following conditions are equivalent:
\begin{enumerate}[label=(\roman*)]
\item $N_a$ is bounded on $H^p$ for at least one $p\in(1,\infty)$;
\item $N_a$ is bounded on $H^2$;
\item $\tau(a)<\infty$;
\item $N_a$ is bounded on $H^p$ for every $p\in(1,\infty)$.
\end{enumerate}
Moreover,
\begin{equation}\label{eq:norm-main}
 \tau(a)\le\|N_a\|_{H^2\to H^2}^2\le C\tau(a)
\end{equation}
with an absolute constant $C$. For each $1<p<\infty$, the $H^p$ norm
is bounded by a finite function of $p$ and $\tau(a)$, and conversely
$\tau(a)$ is bounded by a finite function of $p$ and the $H^p$ norm.
\end{theorem}

Thus boundedness on $H^2$ determines boundedness throughout the Hardy
scale $1<p<\infty$.
For general weights, it remains open whether
$\|N_a\|_{H^p\to H^p}\asymp_p\|N_a\|_{H^2\to H^2}$.

\subsection{The scalar symbol}

The same test can be read from the coefficients of a single analytic
function. Put
\begin{equation}\label{eq:symbol}
 c_n=\frac{a_n}{A_n},\qquad G_a(z)=\sum_{n\ge0}c_nz^n.
\end{equation}
Here $c_0=1$ and $0\le c_n<1$ for $n\ge1$, so $G_a$ is analytic in $\D$.
Define
\begin{equation}\label{eq:truncated-functional}
 \mathcal F_{\mathrm{tr}}(G_a)
 =\sup_{J\ge1}\sum_{k\ge1}
 \min\left\{1,\sum_{j=kJ+1}^{(k+1)J}c_j\right\}^{2}.
\end{equation}
Theorem~\ref{thm:symbol} gives
$\tau(a)\asymp1+\mathcal F_{\mathrm{tr}}(G_a)$ with absolute constants.
Together with Theorem~\ref{thm:main}, this makes
$\mathcal F_{\mathrm{tr}}(G_a)<\infty$ a symbol criterion for every
$H^p$ in the range $1<p<\infty$. The truncation reflects
$0\le R_J(n)\le1$. For doubling partial sums, the block condition
reduces to the familiar $\BMOA$ condition.

\begin{corollary}\label{cor:bmo-intro}
Suppose that, for some $K<\infty$,
\begin{equation}\tag{D}\label{eq:doubling}
 A_{2n+1}\le KA_n\qquad(n\ge0).
\end{equation}
Then $N_a$ is bounded on one, equivalently every, $H^p$, $1<p<\infty$,
if and only if $G_a\in\BMOA$.
Furthermore, without assuming \eqref{eq:doubling},
\[
 G_a\in\BMOA
 \quad\Longleftrightarrow\quad
 \text{$(A_n)$ satisfies \eqref{eq:doubling} and $\tau(a)<\infty$.}
\]
\end{corollary}

Within the bounded family, Corollary~\ref{cor:bloch-symbol} adds
membership of $G_a$ in the Bloch space to these equivalent conditions.
Boundedness of $N_a$ on $\BMOA$ itself is a separate endpoint problem.

\subsection{Proof strategy and notation}

In Section~\ref{sec:l2}, an Abel identity expresses
$\|N_af\|_{H^2}^2$ as the squared norm of a vector-valued Hankel
operator applied to $f$, with an additional squared multiplier norm
when the weights are summable. A reproducing-kernel test for the Hankel
operator then yields \eqref{eq:norm-main}.

To pass from $H^2$ to $H^p$, we need a frequency decomposition adapted
to $(A_n)$. The link between testing and frequency geometry comes from
\eqref{eq:rowsum}. If $A_n\ge2A_{n-J}$, then $R_J(n)\ge1/2$. Thus
at a fixed scale $J$, there are at most $4J\tau(a)$ such indices
$n\ge J$. Section~\ref{sec:low} shows that the resulting counting
condition forces the doubling-time endpoints to have finite lacunary
order. The successor theorem of Sj\"ogren and Sj\"olin
\cite[Theorem~1.2]{SS} and de Leeuw's restriction theorem
\cite[Corollary~2.5]{deLeeuw} turn this geometry into
periodic Littlewood--Paley estimates. In particular, the decomposition
needed in the proof follows from testing.

Section~\ref{sec:hp} uses the square function on the input to prove
sufficiency and a square function on the range of $N_a$ to prove
necessity. Rubio de Francia's inequality handles sufficiency for
$p\le2$ and necessity for $p\ge2$. The adapted Littlewood--Paley
estimates handle the two remaining ranges.
Section~\ref{sec:bmo} proves the symbol criterion and compares it
with Fefferman's positive-coefficient criterion for $\BMOA$.

We use the standard identification of $H^p$ with its boundary functions
in $L^p(\T)$. The operator $S_m$ is the Taylor projection onto degrees
$0,\ldots,m$, and $P_I$ is the Fourier projection onto an interval
$I\subset\mathbb Z$. By the M.~Riesz theorem \cite[Chapter~III, \S2]{Garnett},
$B_p:=\sup_m\|S_m\|_{H^p\to H^p}<\infty$ for $1<p<\infty$.
Set $\Lambda_bf=\sum_nb_n\wh f(n)z^n$ and $\varphi=\sum_na_nz^n$.
The factorization $N_a=\Lambda_{1/A}M_\varphi$ is understood
coefficientwise, with $1/A=(1/A_n)_{n\ge0}$. Each coefficient of
$M_\varphi f$ is a finite sum, so this formula applies to the formal
series $\varphi$. We write $dA$ for normalized area measure and
$X\lesssim_pY$ for $X\le C_pY$, with analogous notation for other
dependencies. Our kernel convention is
$k_w(z)=(1-wz)^{-1}$, so $\|k_w\|_2=(1-|w|^2)^{-1/2}$.

\section{The Hilbert-space testing theorem}\label{sec:l2}

\begin{theorem}[Hilbert-space characterization]\label{thm:l2}
Let $a_n\ge0$, $a_0>0$. The following are equivalent, and
$\|N_a\|_{H^2\to H^2}^2\asymp\tau(a)$ with absolute constants.
\begin{enumerate}
\item[(i)] $N_a$ is bounded on $H^2$.
\item[(ii)] $\tau(a)<\infty$.
\item[(iii)] $\displaystyle\sup_{J\ge1}\frac1J\sum_{n\ge J}\min\Big\{1,\log\frac{A_n}{A_{n-J}}\Big\}^2<\infty$.
\item[(iv)] $\displaystyle\sup_{0<r<1}(1-r)\|N_ak_r\|_{H^2}^2<\infty$.
\end{enumerate}
\end{theorem}

The proof uses an exact Abel identity and a kernel test for vector-valued
Hankel operators. The same representation gives the $H^p$ square
function in Section~\ref{sec:hp}.

We use $H^p(\cE)$ for the Hardy space of Hilbert-space-valued analytic
functions, with norm obtained by taking the $L^p$ norm of the pointwise
$\cE$ norm. The boundary-value and Fourier descriptions agree; see
Blasco \cite{Blasco88,Blasco97} for the vector-valued framework.
In particular, for $F(z)=\sum_{i\ge0}v_i z^i$, Parseval's identity gives
\[
 \|F\|_{H^2(\cE)}^2=\sum_{i\ge0}\|v_i\|_{\cE}^2.
\]
Thus a square-summable $\cE$-valued coefficient sequence is identified
isometrically with its analytic function in $H^2(\cE)$.
Here a vector-valued Hankel operator acts on scalar polynomials by a
matrix $(\bh_{i+k})_{i,k\ge0}$ with entries $\bh_j\in\cE$ and takes
values in $H^2(\cE)$. All the estimates below are independent of the
dimension of $\cE$.

\begin{lemma}[Abel identity and vector-valued Hankel representation]\label{lem:hankel}
Let $a=(a_n)_{n\ge0}$ be a sequence of complex numbers and let
$A_n=\sum_{j=0}^na_j$. Assume that $A_0\ne0$ and that $(|A_n|)_{n\ge0}$
is nondecreasing; in particular, $A_n\ne0$ for every $n$. Put
\[
\delta_m=|A_m|^{-2}-|A_{m+1}|^{-2},\qquad
b=\lim_{n\to\infty}|A_n|^{-2}\ge0.
\]
Thus $\delta_m\ge0$. Write $c=(c_n)_{n\ge0}$, where $c_n=a_n/A_n$,
and initially regard $\varphi(z)=\sum_{n\ge0}a_nz^n$ as a formal power
series.

For a scalar polynomial $f$, set
\[
 y_n=\sum_{k=0}^na_{n-k}\wh f(k),\qquad
 \wh{N_af}(n)=\frac{y_n}{A_n}\quad(n\ge0).
\]
The formal product $\varphi f$ has coefficients $y_n$, and
\[
 S_m(\varphi f)(z):=\sum_{n=0}^my_nz^n
\]
is an ordinary polynomial. If $b>0$, then $(A_n)$ and $(a_n)$ are
bounded, so $\varphi$ is analytic on $\D$ and $\varphi f$ is an actual
analytic product. With the last term omitted when $b=0$, we have
\begin{equation}\label{eq:abel}
 \|N_af\|_{H^2}^2
 =\sum_{m\ge0}\delta_m\|S_m(\varphi f)\|_{H^2}^2
  +b\|\varphi f\|_{H^2}^2.
\end{equation}
All norms and nonnegative sums in this identity are allowed to be
infinite.

If additionally $c\in\ell^2$, the vectors
\[
 \bh_j=(\sqrt{\delta_m}\,a_{m-j})_{m\ge0}\in\ell^2,
 \qquad a_n=0\quad(n<0),
\]
define a vector-valued Hankel operator on scalar polynomials by
\[
 (\Gamma_\bh f)(z)
 =\sum_{i\ge0}\left(\sum_{k\ge0}\bh_{i+k}\wh f(k)\right)z^i
 \in H^2(\ell^2).
\]
Moreover,
\begin{align}
 \|\Gamma_\bh f\|_{H^2(\ell^2)}^2
 &=\sum_{m\ge0}\delta_m\|S_m(\varphi f)\|_{H^2}^2,
       \label{eq:hankel-isometry}\\
 \sum_{j\ge0}\|\bh_j\|_{\ell^2}^2
 &=\sum_{n\ge0}|a_n|^2(|A_n|^{-2}-b)
 \le\|c\|_{\ell^2}^2.
       \label{eq:hankel-symbol}
\end{align}
Consequently, the Abel identity has the norm-decomposition form
\begin{equation}\label{eq:abel-hankel}
 \|N_af\|_{H^2}^2
 =\|\Gamma_\bh f\|_{H^2(\ell^2)}^2
  +b\|\varphi f\|_{H^2}^2,
\end{equation}
again omitting the last term when $b=0$.
In particular, if $|A_n|\to\infty$, then
$\|N_af\|_{H^2}=\|\Gamma_\bh f\|_{H^2(\ell^2)}$.
\end{lemma}

\begin{proof}
For polynomial $f$, the coefficients of $N_af$ are bounded:
indeed, $|a_{n-k}|/|A_n|\le2$ for $0\le k\le n$, and therefore
$|y_n/A_n|\le2\sum_k|\wh f(k)|$. Thus $N_af$ is analytic on $\D$,
although its $H^2$ norm may be infinite.

For $M\ge n$, telescoping gives
\[
 \sum_{m=n}^M\delta_m=|A_n|^{-2}-|A_{M+1}|^{-2}.
\]
Letting $M\to\infty$ yields $|A_n|^{-2}=b+\sum_{m\ge n}\delta_m$.
Consequently, Parseval's identity gives
\begin{align*}
 \|N_af\|_{H^2}^2
 &=\sum_{n\ge0}\frac{|y_n|^2}{|A_n|^2}\\
 &=b\sum_{n\ge0}|y_n|^2
   +\sum_{n\ge0}\sum_{m\ge n}\delta_m|y_n|^2\\
 &=b\sum_{n\ge0}|y_n|^2
   +\sum_{m\ge0}\delta_m\sum_{n=0}^m|y_n|^2.
\end{align*}
Here and below the term involving $b$ is omitted when $b=0$.
The interchange in the last equality is Tonelli's theorem applied to
the nonnegative double series
$\sum_{n\ge0}\sum_{m\ge n}\delta_m|y_n|^2$.
Since
\[
 \sum_{n=0}^m|y_n|^2=\|S_m(\varphi f)\|_{H^2}^2,
 \qquad
 \sum_{n\ge0}|y_n|^2=\|\varphi f\|_{H^2}^2\quad(b>0),
\]
this proves \eqref{eq:abel}, including when its value is infinite.

Now assume $c\in\ell^2$. To establish \eqref{eq:hankel-symbol},
apply Tonelli to the nonnegative double series
$\sum_{j\ge0}\sum_{m\ge j}\delta_m|a_{m-j}|^2$.
Reversing the order of summation, reindexing the finite inner sum by
$n=m-j$, and reversing the order once more gives
\begin{align*}
 \sum_{j\ge0}\|\bh_j\|_{\ell^2}^2
 &=\sum_{j\ge0}\sum_{m\ge j}\delta_m|a_{m-j}|^2\\
 &=\sum_{m\ge0}\delta_m\sum_{j=0}^m|a_{m-j}|^2
  =\sum_{m\ge0}\delta_m\sum_{n=0}^m|a_n|^2\\
 &=\sum_{n\ge0}|a_n|^2\sum_{m\ge n}\delta_m\\
 &=\sum_{n\ge0}|a_n|^2\bigl(|A_n|^{-2}-b\bigr)
 \le\sum_{n\ge0}\left|\frac{a_n}{A_n}\right|^2
 =\|c\|_{\ell^2}^2<\infty.
\end{align*}
Thus each $\bh_j$ belongs to $\ell^2$.
For each $k$, the coefficient sequence $(\bh_{i+k})_{i\ge0}$ is a
square-summable tail of $(\bh_j)_{j\ge0}$ and so defines a function in
$H^2(\ell^2)$. Finite linear combinations of these functions define
$\Gamma_\bh f$ for every scalar polynomial $f$.

To verify \eqref{eq:hankel-isometry}, extend $y_n$ by zero for $n<0$.
The $m$th scalar coordinate of the $i$th Taylor coefficient of
$\Gamma_\bh f$ is
\[
 \bigl[\wh{\Gamma_\bh f}(i)\bigr]_m
 =\sqrt{\delta_m}\sum_{k\ge0}a_{m-i-k}\wh f(k)
 =\sqrt{\delta_m}\,y_{m-i}.
\]
It vanishes when $m<i$. Parseval's identity in $H^2(\ell^2)$ and
Tonelli's theorem now give
\begin{align*}
 \|\Gamma_\bh f\|_{H^2(\ell^2)}^2
 &=\sum_{i\ge0}\sum_{m\ge i}\delta_m|y_{m-i}|^2\\
 &=\sum_{m\ge0}\delta_m\sum_{i=0}^m|y_{m-i}|^2\\
 &=\sum_{m\ge0}\delta_m\sum_{n=0}^m|y_n|^2
  =\sum_{m\ge0}\delta_m\|S_m(\varphi f)\|_{H^2}^2.
\end{align*}
Combining this with \eqref{eq:abel} proves \eqref{eq:abel-hankel}.
If $|A_n|\to\infty$, then $b=0$, giving the final assertion.
\end{proof}

\begin{remark}
The hypothesis $c\in\ell^2$ cannot be suppressed in the Hilbert-space
formulation. For $a_n=2^n$, one has
$\delta_m|a_m|^2\to3/16$, so $\bh_0\notin\ell^2$.
The extended-sum identity \eqref{eq:abel} remains valid.
\end{remark}

The next lemma is a vector-valued version of Bonsall's scalar Hankel
kernel test \cite{Bonsall}. Its proof uses the Littlewood--Paley and
Carleson-measure characterizations of BMOA; see
\cite[Chapter~VI, \S3]{Garnett} and \cite[\S\S1--2]{Blasco97}
for the scalar and Hilbert-valued formulations, respectively.
We include the proof to specify the dimension-independent reduction.

\begin{lemma}[Kernel test for vector-valued Hankel operators]\label{lem:rkt}
Let $\cE$ be a complex Hilbert space and let
$(\bh_j)_{j\ge0}$ be a sequence in $\cE$ such that
$\sum_{j\ge0}\|\bh_j\|_{\cE}^2<\infty$. For a scalar polynomial $f$,
define the vector-valued Hankel operator by
\[
 (\Gamma_\bh f)(z)
 =\sum_{i\ge0}\left(\sum_{k\ge0}\bh_{i+k}\wh f(k)\right)z^i.
\]
This is a function in $H^2(\cE)$. For $w\in\D$, its value on the kernel
$k_w(z)=(1-wz)^{-1}$ is defined by the norm-convergent series
\[
 \Gamma_\bh k_w
 :=\sum_{k\ge0}w^k\left(\sum_{i\ge0}\bh_{i+k}z^i\right)
 \quad\text{in }H^2(\cE).
\]
Put
\[
 K=\sup_{w\in\D}(1-|w|^2)^{1/2}
       \|\Gamma_\bh k_w\|_{H^2(\cE)}.
\]
Then $\Gamma_\bh$ extends to a bounded operator from $H^2$ into
$H^2(\cE)$ if and only if $K<\infty$. In that case,
\[
 K\le\|\Gamma_\bh\|_{H^2\to H^2(\cE)}\le CK,
\]
where $C$ is an absolute constant, independent of $\cE$.
\end{lemma}

\begin{proof}
Write $H=(\sum_{j\ge0}\|\bh_j\|_{\cE}^2)^{1/2}$.
For each $k\ge0$, denote the image of the monomial $z^k$ by
\[
 \mathcal C_k(z)=\sum_{i\ge0}\bh_{i+k}z^i,
 \qquad
 \|\mathcal C_k\|_{H^2(\cE)}^2
 =\sum_{i\ge0}\|\bh_{i+k}\|_{\cE}^2\le H^2.
\]
Thus $\Gamma_\bh f=\sum_k\wh f(k)\mathcal C_k$ is well defined for
every scalar polynomial $f$, and
\[
 \sum_{k\ge0}|w|^k\|\mathcal C_k\|_{H^2(\cE)}
 \le\frac{H}{1-|w|}<\infty.
\]
This proves norm convergence of the series defining $\Gamma_\bh k_w$
without assuming boundedness of $\Gamma_\bh$.
If a bounded extension exists, this series agrees with its value on
$k_w$, since the Taylor polynomials of $k_w$ converge in $H^2$.
As $\|k_w\|_{H^2}=(1-|w|^2)^{-1/2}$, testing on normalized kernels gives
$K\le\|\Gamma_\bh\|$.

\smallskip
\noindent\emph{Poisson oscillation.}
Assume $K<\infty$. All outputs lie in the separable
closed span of the vectors $\bh_j$, so we may replace $\cE$ by this
span. Define
\[
 G(z)=\sum_{j\ge0}\bh_jz^{j+1}.
\]
Then $G\in H^2(\cE)$, $G(0)=0$, and $\|G\|_{H^2(\cE)}=H$.
The elementary identity
\[
 \frac{z^{j+1}-w^{j+1}}{z-w}=\sum_{i=0}^jz^iw^{j-i}
\]
gives, for $z,w\in\D$ (with the removable value at $z=w$),
\[
 \frac{G(z)-G(w)}{z-w}
 =\sum_{i\ge0}\left(\sum_{k\ge0}\bh_{i+k}w^k\right)z^i.
\]
The function on the right is precisely $\Gamma_\bh k_w$.
Parseval's identity and $|\zeta-w|=|1-\bar w\zeta|$ for $\zeta\in\T$
therefore yield
\begin{equation}\label{eq:garsia}
 (1-|w|^2)\|\Gamma_\bh k_w\|_{H^2(\cE)}^2
 =\int_\T\|G(\zeta)-G(w)\|_{\cE}^2P_w(\zeta)\,dm(\zeta),
\end{equation}
where $P_w(\zeta)=(1-|w|^2)/|1-\bar w\zeta|^2$ is the Poisson kernel.
In particular, the integral in \eqref{eq:garsia} is at most $K^2$.

\smallskip
\noindent\emph{Carleson embedding.}
The scalar Littlewood--Paley identity
\cite[Chapter~VI, Lemma~3.1]{Garnett} reads
\[
 \|F\|_{H^2}^2-|F(0)|^2
 =\int_\D|F'(z)|^2\log\frac1{|z|^2}\,dA(z).
\]
For an $\cE$-valued function $F$, apply this identity to its scalar
coordinates $F_\ell$ in an orthonormal basis. Tonelli's theorem gives
\[
 \sum_\ell\int_\D|F_\ell'(z)|^2\log\frac1{|z|^2}\,dA(z)
 =\int_\D\|F'(z)\|_{\cE}^2\log\frac1{|z|^2}\,dA(z).
\]
Summing the left sides of the scalar identities likewise gives
$\|F\|_{H^2(\cE)}^2-\|F(0)\|_{\cE}^2$. Hence
\begin{equation}\label{eq:rkt-lp}
 \|F\|_{H^2(\cE)}^2-\|F(0)\|_{\cE}^2
 =\int_\D\|F'(z)\|_{\cE}^2\log\frac1{|z|^2}\,dA(z),
\end{equation}
with the same constant as in the scalar case.

Let $\sigma_w(z)=(w-z)/(1-\bar wz)$.
The boundary change of variables under $\sigma_w$ identifies the
integral in \eqref{eq:garsia} with
$\|G\circ\sigma_w-G(w)\|_{H^2(\cE)}^2$.
Apply \eqref{eq:rkt-lp} to $G\circ\sigma_w$ and then change variables
using $\sigma_w^{-1}=\sigma_w$. The derivative factor is canceled by
the area Jacobian, giving
\begin{equation}\label{eq:rkt-green}
 J(w):=\int_\D\|G'(z)\|_{\cE}^2
             \log\frac1{|\sigma_w(z)|^2}\,dA(z)
 \le K^2.
\end{equation}

Define the positive scalar measure
\[
 d\nu(z)=\|G'(z)\|_{\cE}^2L(z)\,dA(z),
 \qquad L(z)=\log\frac1{|z|^2}.
\]
The case $w=0$ in \eqref{eq:rkt-green} gives $\nu(\D)\le K^2$.
We claim that $\nu$ is a Carleson measure with constant at most $CK^2$.
For an arc $I\subset\T$, let $|I|=m(I)$ and define
\[
 S(I)=\{r\zeta:\zeta\in I,\ 1-|I|\le r<1\}.
\]
Suppose first that $0<|I|\le1/2$, and choose
$w=(1-|I|)\zeta_I$, where $\zeta_I$ is the midpoint of $I$.
For $z\in S(I)$, elementary geometry gives
\[
 |1-\bar wz|\le C|I|,\qquad 1-|w|^2\ge |I|.
\]
Using $\log(1/t)\ge1-t$ for $0<t\le1$ and the automorphism identity,
we obtain
\[
 \log\frac1{|\sigma_w(z)|^2}
 \ge 1-|\sigma_w(z)|^2
 =\frac{(1-|w|^2)(1-|z|^2)}{|1-\bar wz|^2}
 \ge c\,\frac{1-|z|^2}{|I|}.
\]
On $S(I)$ we have $|z|\ge1/2$, so $L(z)\asymp1-|z|^2$.
Integrating the preceding inequality over $S(I)$ and using
\eqref{eq:rkt-green} gives
\[
 \nu(S(I))\le C|I|J(w)\le CK^2|I|.
\]
For $|I|>1/2$, the same bound follows from $\nu(\D)\le K^2$.
This proves the claim.

The scalar Carleson embedding theorem
\cite[Theorem~9.3]{Duren} now implies
\begin{equation}\label{eq:rkt-embedding}
 \int_\D|q(z)|^2\,d\nu(z)\le CK^2\|q\|_{H^2}^2
 \qquad(q\in H^2).
\end{equation}
Its Hilbert-valued version follows directly by applying this scalar
estimate to the coordinates of $Y=\sum_\ell Y_\ell e_\ell$:
\[
 \int_\D\|Y(z)\|_{\cE}^2\,d\nu(z)
 =\sum_\ell\int_\D|Y_\ell(z)|^2\,d\nu(z)
 \le CK^2\sum_\ell\|Y_\ell\|_{H^2}^2
 =CK^2\|Y\|_{H^2(\cE)}^2.
\]
Here Tonelli justifies the equality, and the constant is independent
of the number of coordinates.

\smallskip
\noindent\emph{The Hankel pairing.}
Take a scalar polynomial $f(z)=\sum_kx_kz^k$ and an $\cE$-valued
polynomial $Y(z)=\sum_i\by_i z^i$, and put
\[
 q(z)=\sum_k\overline{x_k}z^k,\qquad \Phi(z)=zq(z)Y(z).
\]
Thus $x_k$ and $\by_i$ are the Taylor coefficients of the functions
$f$ and $Y$, respectively, and $\|q\|_{H^2}=\|f\|_{H^2}$.
Take inner products linear in the first variable.
The coefficient of $z^{j+1}$ in $\Phi$ is
$\sum_{i+k=j}\overline{x_k}\by_i$, so
\[
 \langle\Gamma_\bh f,Y\rangle_{H^2(\cE)}
 =\sum_{i,k}x_k\langle\bh_{i+k},\by_i\rangle_{\cE}
 =\langle G,\Phi\rangle_{H^2(\cE)}.
\]
Polarizing \eqref{eq:rkt-lp}, and using $\Phi(0)=0$, gives
\[
 \langle G,\Phi\rangle_{H^2(\cE)}
 =\int_\D\langle G'(z),\Phi'(z)\rangle_{\cE}L(z)\,dA(z).
\]
Since $\Phi'=(zq)'Y+zqY'$, Cauchy--Schwarz applied separately to these
two terms yields
\begin{align*}
 |\langle\Gamma_\bh f,Y\rangle|
 &\le
 \left(\int_\D\|Y\|_{\cE}^2\,d\nu\right)^{1/2}
 \left(\int_\D|(zq)'|^2L\,dA\right)^{1/2}\\
 &\quad+
 \left(\int_\D|zq|^2\,d\nu\right)^{1/2}
 \left(\int_\D\|Y'\|_{\cE}^2L\,dA\right)^{1/2}.
\end{align*}
The scalar and Hilbert-valued Carleson embeddings control the first
factor in each product. The scalar Littlewood--Paley identity and
\eqref{eq:rkt-lp} control the second factors:
\[
 \int_\D|(zq)'|^2L\,dA=\|zq\|_{H^2}^2=\|q\|_{H^2}^2,
 \qquad
 \int_\D\|Y'\|_{\cE}^2L\,dA\le\|Y\|_{H^2(\cE)}^2.
\]
Consequently,
\[
 |\langle\Gamma_\bh f,Y\rangle|
 \le CK\|q\|_{H^2}\|Y\|_{H^2(\cE)}
 =CK\|f\|_{H^2}\|Y\|_{H^2(\cE)}.
\]
Taking the supremum over $\cE$-valued polynomials $Y$ of
$H^2(\cE)$ norm at most one gives
$\|\Gamma_\bh f\|_{H^2(\cE)}\le CK\|f\|_{H^2}$.
Density of the scalar polynomials then supplies the unique bounded
extension from $H^2$ into $H^2(\cE)$.
\end{proof}

\begin{proof}[Proof of Theorem~\ref{thm:l2}]
Testing on $D_J$, for which $\|D_J\|_{H^2}^2=J$, proves
(i)$\Rightarrow$(ii), with $\tau(a)\le\|N_a\|_{H^2\to H^2}^2$. For (ii)$\Leftrightarrow$(iii), use
\[
(1-e^{-1})\min(1,t)\le1-e^{-t}\le\min(1,t),
\qquad t=\log(A_n/A_{n-J}).
\]

Write
\[
 \tau'=\sup_{0<r<1}(1-r)\|N_ak_r\|_{H^2}^2.
\]
We prove the quantitative comparison
\[
 \tau'\le\tau(a)\le e^2\tau',
\]
which gives (ii)$\Leftrightarrow$(iv).

First assume $\tau(a)<\infty$ and fix $0<r<1$. The coefficient of
$z^n$ in $D_J$ is one for $J\ge n+1$ and zero otherwise. Since
$(1-r)\sum_{J\ge n+1}r^{J-1}=r^n$, we have the expansion
\[
 k_r(z)=\frac1{1-rz}
       =(1-r)\sum_{J\ge1}r^{J-1}D_J(z).
\]
The weights $(1-r)r^{J-1}$ sum to one. Since each output coefficient
in \eqref{eq:operator} uses only finitely many input coefficients,
$\wh{N_ak_r}(n)$ is the same weighted average of $\wh{N_aD_J}(n)$.
Apply Jensen's inequality for the squared modulus to that average
and sum over $n$. Exchanging the sums over $n$ and $J$ by Tonelli's
theorem gives
\begin{align*}
 \|N_ak_r\|_{H^2}^2
 &\le(1-r)\sum_{J\ge1}r^{J-1}\|N_aD_J\|_{H^2}^2\\[4pt]
 &\le\tau(a)(1-r)\sum_{J\ge1}Jr^{J-1}
 =\frac{\tau(a)}{1-r}.
\end{align*}
Taking the supremum over $r$ proves $\tau'\le\tau(a)$.

Conversely, assume $\tau'<\infty$. Fix $J\ge2$ and set $r_J=1-1/J$.
For $0\le n<J$,
\[
 r_J^n\ge r_J^{J-1}
 =\left(1+\frac1{J-1}\right)^{-(J-1)}
 \ge e^{-1},
\]
where the last inequality follows from $\log(1+t)\le t$.
Every Taylor coefficient of $D_J$ is therefore at most $e$ times
the corresponding coefficient of $k_{r_J}$.
Positivity of $N_a$ preserves this coefficientwise inequality, so
\[
 \frac{\|N_aD_J\|_{H^2}^2}{J}
 \le\frac{e^2}{J}\|N_ak_{r_J}\|_{H^2}^2
 =e^2(1-r_J)\|N_ak_{r_J}\|_{H^2}^2
 \le e^2\tau'.
\]
For $J=1$, the coefficients of $D_1=1$ are bounded by those of $k_r$
for every $0<r<1$. Positivity again gives
$\|N_aD_1\|_{H^2}^2\le\tau'/(1-r)$.
Letting $r\downarrow0$ gives $\|N_aD_1\|_{H^2}^2\le\tau'$, and hence
$\tau(a)\le e^2\tau'$.

It remains to prove boundedness from (iv). Since
$\wh{N_aD_1}(n)=c_n$, the estimate for $J=1$ also shows that
$c=(c_n)_{n\ge0}\in\ell^2$.
For every $w\in\D$ and every $n\ge0$, positivity gives
\[
 |\wh{N_ak_w}(n)|\le\wh{N_ak_{|w|}}(n).
\]
Consequently,
\[
 \sup_{w\in\D}(1-|w|^2)\|N_ak_w\|_{H^2}^2\le2\tau'.
\]
The coefficient calculation in Lemma~\ref{lem:hankel} applies also
to $k_w$: each convolution coefficient is a finite sum, and the
subsequent norm identities follow by Tonelli. Thus
\[
 \|\Gamma_\bh k_w\|_{H^2(\ell^2)}
 \le\|N_ak_w\|_{H^2}.
\]
Lemma~\ref{lem:rkt}, with $\cE=\ell^2$, therefore gives
\[
 \|\Gamma_\bh\|_{H^2\to H^2(\ell^2)}^2\le2C^2\tau'.
\]

We now use the norm decomposition \eqref{eq:abel-hankel}.
If $A_n\to\infty$, then $b=0$, and for every scalar polynomial $f$,
\[
 \|N_af\|_{H^2}^2
 =\|\Gamma_\bh f\|_{H^2(\ell^2)}^2
 \le2C^2\tau'\|f\|_{H^2}^2.
\]
If instead $A_n\to A_\infty<\infty$, nonnegativity of the weights gives
$\varphi\in H^\infty$ and $\|\varphi\|_\infty=A_\infty$.
Indeed, $|\varphi(z)|\le\sum_na_n=A_\infty$, while
$\varphi(r)\uparrow A_\infty$ as $r\uparrow1$.
The multiplication operator
\[
 (M_{\varphi/A_\infty}f)(z)
 =\frac{\varphi(z)}{A_\infty}f(z)
\]
is therefore a contraction on $H^2$.
Since $b=A_\infty^{-2}$, \eqref{eq:abel-hankel} becomes
\[
 \|N_af\|_{H^2}^2
 =\|\Gamma_\bh f\|_{H^2(\ell^2)}^2
  +\|M_{\varphi/A_\infty}f\|_{H^2}^2
 \le(2C^2\tau'+1)\|f\|_{H^2}^2.
\]
In either case, density of the polynomials gives a bounded operator
$N_a$ on $H^2$. Since $\tau'\le\tau(a)$ and $\tau(a)\ge1$, these
estimates also prove the required absolute upper bound for its norm.
\end{proof}

\begin{remark}[Relation to the discrete Weiss conjecture]
There is also an admissibility interpretation of the kernel test.
The weighted shift $V_ae_n=(A_n/A_{n+1})e_{n+1}$ is a contraction, and
$N_ae_k=V_a^kc$, where $c=(a_n/A_n)_{n\ge0}\in\ell^2$ whenever the
test is finite. Hence boundedness of $N_a$ is equivalent to
\[
 \sum_{k\ge0}|\langle V_a^{*k}y,c\rangle|^2\lesssim\|y\|_2^2.
\]
Moreover $N_ak_w=(I-wV_a)^{-1}c$. The equivalence between the
normalized kernel bound and boundedness is therefore an instance of
the discrete Weiss conjecture for contractions with scalar output;
see Harper \cite[Theorem~4.4]{Harper}. Positivity reduces the test to real $w$.
The preceding proof also supplies the explicit vector-valued Hankel
representation used below.
\end{remark}

\section{Geometry of the doubling-time partition}\label{sec:low}

Throughout this section assume $A_n\to\infty$. Define the doubling times
and their frequency blocks by
\[
 n_0=0,\qquad n_{j+1}=\min\{n:A_n\ge2A_{n_j}\},\qquad
 W_j=[n_j,n_{j+1})\cap\mathbb Z.
\]
Each block is finite, although its length may be arbitrarily large.
For $1<q<\infty$, the property needed below is
\begin{equation}\tag{$\mathrm{LP}_q$}\label{LP}
 \left\|\left(\sum_j|P_{W_j}g|^2\right)^{1/2}\right\|_{L^q(\T)}
 \le L\|g\|_{L^q(\T)}
\end{equation}
for every trigonometric polynomial $g$. We derive this inequality by
using a necessary counting condition for boundedness of $N_a$ to control
the lacunary order of the endpoints $\{n_j\}$.

\subsection{A necessary counting condition}
Recall $R_J$ from \eqref{eq:rowsum}, and put
\begin{equation}\label{eq:count}
 \kappa_0(a)=\sup_{J\ge1}\frac1J
       \#\{n\ge J:A_n\ge2A_{n-J}\}.
\end{equation}

\begin{proposition}\label{prop:nec}
One has $\kappa_0(a)\le4\tau(a)$. If $N_a$ is bounded on $H^p$ for
some $1<p\le2$, then
\[
 \sup_{J\ge1}\frac1J\sum_{n\ge J}R_J(n)^{p'}
 \le C_p\|N_a\|_{H^p\to H^p}^{p'},\qquad
 \kappa_0(a)\le C_p'\|N_a\|_{H^p\to H^p}^{p'},
\]
where $p'=p/(p-1)$.
\end{proposition}
\begin{proof}
On the set counted in \eqref{eq:count}, $R_J(n)\ge1/2$; this proves
the first assertion. For the second, set $F=N_aD_J$. The standard
Dirichlet-kernel estimate $\|D_J\|_p\asymp_pJ^{1/p'}$ and the
Hausdorff--Young inequality (see \cite[Chapter~3]{Grafakos}) give
\[
 \sum_{n\ge0}|\wh F(n)|^{p'}\le\|F\|_p^{p'}
 \le C_p\|N_a\|_{H^p}^{p'}J.
\]
Apply \eqref{eq:rowsum} and again use $R_J\ge1/2$ on the counted set.
\end{proof}

\subsection{Clustering and finite lacunary order}
Let $\mathcal W=\{W_j\}$ be a partition of $\mathbb N_0$ into consecutive
intervals with left endpoints $Z=\{n_j\}$ and lengths $\ell_j$. A last
block may have infinite length. For an integer $J\ge1$, a block is
\emph{$J$-short} if $\ell_j\le J$. A \emph{$J$-component} is the union
of a maximal run of consecutive $J$-short blocks. We consider the property
\begin{equation}\tag{$\mathrm P_M$}\label{PM}
 \text{At each integer scale }J\ge1\text{ there are at most }M\
 J\text{-components, each of length at most }MJ.
\end{equation}

\begin{lemma}[Clustering]\label{lem:cluster}
If $\kappa_0(a)<\infty$, the doubling-time partition satisfies
$(\mathrm P_M)$ with $M=\lceil2\kappa_0(a)+2\rceil$.
\end{lemma}

\begin{proof}
Fix $J\ge1$. If $W_i$ is $J$-short and $n_i\ge J$, then every integer
$n\in[n_i+J,n_i+2J)$ satisfies $n\ge2J$, $n\ge n_{i+1}$, and
$n-2J\le n_i$. Hence
\[
 A_n\ge A_{n_{i+1}}\ge2A_{n_i}\ge2A_{n-2J}.
\]
By \eqref{eq:count}, the union of these intervals contains at most
$2J\kappa_0(a)$ integers.

Each translated block $W_i+J$ is contained in its corresponding
interval $[n_i+J,n_i+2J)$. The translated blocks are disjoint, so the
total length of the $J$-short blocks starting at $n_i\ge J$ is at most
$2J\kappa_0(a)$.

Choose one such block from each $J$-component that contains one. Any
two chosen blocks are separated by a block of length greater than $J$.
Their starts therefore differ by more than $J$, so their intervals
$[n_i+J,n_i+2J)$ are disjoint and contribute $J$ integers each. This
bounds the number of these components by $2\kappa_0(a)$.

The $J$-short blocks starting before $J$ lie in $[0,2J)$ and belong to
at most one component, since starts in different components differ by
more than $J$. Thus there are at most $2\kappa_0(a)+1$ components, each of
length at most $(2\kappa_0(a)+2)J$.
\end{proof}

\begin{proposition}[Finite lacunary order]\label{prop:order}
Let $M\ge1$ be an integer. If $\mathcal W$ satisfies $(\mathrm P_M)$, there are sets $E_0\subset E_1\subset\dots\subset E_M=Z$, with
$E_0$ a single point, such that for $1\le i\le M$
\begin{equation}\label{eq:succ}
|x-y|\ge\frac1M\operatorname{dist}(x,E_{i-1})\qquad\text{for all }x,y\in E_i,\ x\ne y .
\end{equation}
In the terminology of Sj\"ogren and Sj\"olin \cite{SS}, each $E_i$ is a successor of $E_{i-1}$
(with constant $1/2M$), and $Z$ is a lacunary set of order at most $M$.
\end{proposition}

\begin{proof}
If $Z$ is a singleton, take $E_i=Z$ for every $i$. Otherwise, assign
parents to endpoints as their components merge. All strict ancestors
of an endpoint will be counted at the first scale at which the endpoint
belongs to a nontrivial component.

For $x\in Z$, let
\[
 s(x)=\operatorname{dist}(x,Z\setminus\{x\}).
\]
This is a positive integer. Order $Z$ lexicographically by $(s(x),x)$;
this is a well-order, denoted by $\prec$, and we write $x_*$ for its
least element.

For an integer $J\ge1$, join consecutive endpoints whenever the block
between them has length at most $J$. Let $C_J(x)$ consist of all endpoints
reachable from $x$ along these edges, including $x$ itself. Choose its
representative by setting
\[
 r_J(x)=\min\nolimits_\prec C_J(x).
\]
Each nontrivial $C_J(x)$ is the endpoint set of a $J$-component of
$\mathcal W$, with diameter equal to the length of that $J$-component.
By $(\mathrm P_M)$ there are at most $M$ nontrivial endpoint components
at every scale, and
\[
 \operatorname{diam}C_J(x)\le MJ.
\]
For $J<s(x)$ the component $C_J(x)$ is a singleton; for $J\ge s(x)$
it is nontrivial.

For $x\ne x_*$, define its change scale and its parent by
\[
 d(x)=\min\{J\ge s(x):r_J(x)\ne x\},
 \qquad p(x)=r_{d(x)}(x),
\]
and set $d(x_*)=\infty$.
The scale $d(x)$ is finite, since $x$ and $x_*$ lie in a common
component for all sufficiently large $J$.
As the components increase with $J$, a point that ceases to be their
least element never becomes their least element again. Hence
\[
 p(x)\prec x,\qquad
 s(p(x))\le s(x),\qquad
 d(p(x))>d(x).
\]
The first two assertions follow from the definition of $\prec$.
For the third, $p(x)$ is still the representative of its component
at scale $d(x)$.
Following parents gives a strictly decreasing $\prec$-chain, which
must terminate at $x_*$. Let $h(x)$ denote its length, so that
$h(x_*)=0$ and $h(p(x))=h(x)-1$.

\smallskip
To bound $h(x)$, take $J=s(x)$ for $x\ne x_*$. Since $s$ is nonincreasing
and $d$ is strictly increasing along parent links, every strict ancestor
$v$ of $x$ satisfies
\[
 s(v)\le J=s(x)\le d(x)<d(v).
\]
The first inequality places $v$ in a nontrivial component at scale $J$.
The last inequality ensures that $v$ still represents that component.
Distinct ancestors therefore represent distinct nontrivial components
at the same scale. There are at most $M$ such components, so $h(x)\le M$.

We also record the compatibility between representatives and parents:
for $J\ge s(y)$, the point $r_J(y)$ is an ancestor of $y$, possibly
$y$ itself. Indeed, whenever $d(y)\le J$, the points $y$ and $p(y)$
belong to the same $J$-component. Following parents until reaching
a point $v$ with $d(v)>J$ therefore stays in that component and ends
at its representative.

Now set
\[
 E_i=\{x\in Z:h(x)\le i\},\qquad 0\le i\le M.
\]
Then $E_0=\{x_*\}$, the sets increase, and $E_M=Z$.
To prove \eqref{eq:succ}, it suffices to consider
$x\in E_i\setminus E_{i-1}$ and $y\in E_i\setminus\{x\}$.
Since $p(x)\in E_{i-1}$,
\[
 \operatorname{dist}(x,E_{i-1})
 \le |x-p(x)|
 \le \operatorname{diam}C_{d(x)}(x)
 \le Md(x).
\]
We show that $|x-y|\ge d(x)$.
Otherwise the integer $\rho=|x-y|$ satisfies
\[
 \max\{s(x),s(y)\}\le\rho<d(x).
\]
All gaps between $x$ and $y$ have length at most $\rho$, so
$C_\rho(x)=C_\rho(y)$. Since $\rho<d(x)$, their representative is $x$:
\[
 r_\rho(y)=r_\rho(x)=x.
\]
The compatibility observation makes $x$ a strict ancestor of $y$,
giving $h(y)>h(x)=i$, a contradiction.
Consequently,
\[
 \operatorname{dist}(x,E_{i-1})\le Md(x)\le M|x-y|,
\]
as required. Taking the smaller constant $1/(2M)$ gives the strict
successor inequality in the terminology of \cite{SS}.
\end{proof}

\subsection{The Littlewood--Paley estimate}

\begin{theorem}[Littlewood--Paley property]\label{thm:LP}
Suppose $A_n\to\infty$ and $\kappa_0(a)<\infty$. Then the
doubling-time partition $\{W_j\}$ satisfies $(\mathrm{LP}_q)$ for
every $1<q<\infty$, with a constant $L=L(q,\kappa_0(a))$.
\end{theorem}

\begin{proof}
Fix $1<q<\infty$ and put $M=\lceil2\kappa_0(a)+2\rceil$.
Lemma~\ref{lem:cluster} and Proposition~\ref{prop:order} give a
successor chain ending at $Z$, of length at most $M$ and with
successor constants $1/(2M)$. Every set in the chain is closed and
null, being a subset of the integers. Translating the chain by
$-1/2$ preserves these properties and constants. Hence
\cite[Theorem~1.2]{SS} gives the real-line Littlewood--Paley property
for
\[
 E=Z-\tfrac12.
\]

The complementary intervals of $E$ are
\[
 I_-=(-\infty,-\tfrac12),\qquad
 I_j=(n_j-\tfrac12,n_{j+1}-\tfrac12)\quad(j\ge0).
\]
Write $\mathcal I=\{I_-,I_0,I_1,\ldots\}$ and let $S_I$ denote the
real-line Fourier projection with symbol $\mathbf1_I$.
The preceding Littlewood--Paley property means that, for some
$B=B(q,M)\ge1$,
\[
 B^{-1}\|F\|_{L^q(\mathbb R)}
 \le\|Q_{\mathbb R}F\|_{L^q(\mathbb R)}
 \le B\|F\|_{L^q(\mathbb R)},
 \qquad
 Q_{\mathbb R}F
 :=\left(\sum_{I\in\mathcal I}|S_IF|^2\right)^{1/2}.
\]

We first deduce a uniform estimate for signed sums of projections.
For a finite set of block indices $\mathcal J$ and signs
$\epsilon_j\in\{-1,1\}$, put
\[
 T_\epsilon F=\sum_{j\in\mathcal J}\epsilon_jS_{I_j}F.
\]
This is a finite sum of bounded interval projections.
Since the intervals are disjoint,
\[
 Q_{\mathbb R}(T_\epsilon F)
 =\left(\sum_{j\in\mathcal J}|S_{I_j}F|^2\right)^{1/2}
 \le Q_{\mathbb R}F.
\]
Using both sides of the real-line estimate therefore gives
\[
 \|T_\epsilon F\|_{L^q(\mathbb R)}
 \le B\|Q_{\mathbb R}(T_\epsilon F)\|_{L^q(\mathbb R)}
 \le B^2\|F\|_{L^q(\mathbb R)}.
\]
The bound is independent of $\mathcal J$ and the signs.

We now transfer this estimate to $\mathbb T$.
The symbol of $T_\epsilon$ is
\[
 m_\epsilon(\xi)=\sum_{j\in\mathcal J}\epsilon_j\mathbf1_{I_j}(\xi).
\]
Its possible jumps occur only at half-integers, so it is continuous
at every integer. Moreover,
\[
 I_j\cap\mathbb Z=[n_j,n_{j+1})\cap\mathbb Z=W_j,
 \qquad
 m_\epsilon(n)
 =\sum_{j\in\mathcal J}\epsilon_j\mathbf1_{W_j}(n).
\]
At each jump, assign $m_\epsilon$ the average of its two one-sided
limits. This preserves the real-line operator and all integer values,
and makes the symbol regulated in de Leeuw's sense: its symmetric
interval averages converge to its value at every point.
By de Leeuw's restriction theorem \cite[Corollary~2.5, p.~372]{deLeeuw},
the same symbol defines a multiplier of equal norm on
$L^q(\widehat{\mathbb R_d})$, where $\mathbb R_d$ is the additive
group $\mathbb R$ with the discrete topology and $\widehat{\mathbb R_d}$
is its compact dual. For trigonometric polynomials with integer
frequencies, the Bohr mean of the $q$th power of the modulus is the
normalized mean over one period. Thus restricting the symbol to
$\mathbb Z$ gives
\[
 \left\|\sum_{j\in\mathcal J}\epsilon_jP_{W_j}g\right\|_{L^q(\mathbb T)}
 \le B^2\|g\|_{L^q(\mathbb T)}
\]
for every trigonometric polynomial $g$.

Finally, take the signs to be independent Rademacher variables.
Khintchine's inequality \cite[Appendix~C, \S C.2]{Grafakos}, followed by
integration over $\mathbb T$, gives
\begin{align*}
 \left\|\left(\sum_{j\in\mathcal J}|P_{W_j}g|^2\right)^{1/2}
       \right\|_{L^q(\mathbb T)}^q
 &\lesssim_q
 \mathbb E\left\|\sum_{j\in\mathcal J}\epsilon_jP_{W_j}g
              \right\|_{L^q(\mathbb T)}^q\\
 &\le B^{2q}\|g\|_{L^q(\mathbb T)}^q.
\end{align*}
Only finitely many blocks meet the Fourier support of $g$.
Choosing $\mathcal J$ to contain those blocks proves
$(\mathrm{LP}_q)$ with $L\le C_qB(q,M)^2$, which depends only on
$q$ and $\kappa_0(a)$.
\end{proof}

If $A_\infty<\infty$, there are only finitely many doubling times,
say $Z=\{n_0,\ldots,n_N\}$. With the convention $n_{N+1}=+\infty$,
the same interval formulas give the last block
$W_N=[n_N,\infty)\cap\mathbb Z$ and the last complementary interval
$I_N=(n_N-\tfrac12,\infty)$, so $I_N\cap\mathbb Z=W_N$.
The proof of Theorem~\ref{thm:main} treats finite total mass by
regularization.

\section{Square functions and the Hardy-space characterization}\label{sec:hp}

Except in the final regularization argument, assume $A_n\to\infty$.
For an analytic polynomial $f$, define the square function
\[
 \cS f:=\left(\sum_{m\ge0}\delta_m|S_m(\varphi f)|^2\right)^{1/2},
 \qquad
 \delta_m=\frac1{A_m^2}-\frac1{A_{m+1}^2}.
\]
We prove that boundedness of $N_a$ on $H^2$ is equivalent to boundedness
on $H^p$ for $1<p<\infty$. The Hankel representation connects $\cS$
with $H^2$-boundedness. Sufficiency estimates $N_af$ in terms of
$\cS f$, and necessity estimates $\cS f$ in terms of $N_af$.

Write $p'=p/(p-1)$ for the H\"older conjugate of $p$, so
$1/p+1/p'=1$.
At $p=2$, \eqref{eq:abel} gives $\|\cS f\|_2=\|N_af\|_2$.
Theorem~\ref{thm:LP} supplies the partition estimates from the
counting condition in the completion of the proof.

\subsection{Fourier estimates}
We collect the following standard Fourier estimates in the forms needed
below. They follow from the M.~Riesz theorem, Rubio de Francia's
inequality, and the vector-valued Carleson theorem; the proof records
the reductions and references. All constants are independent of the
number and location of the intervals.

\begin{lemma}[Fourier estimates]\label{lem:fourier}
Let $1<q<\infty$. All scalar functions below are complex-valued
functions on $\T$, and $\|\cdot\|_q$ denotes the $L^q(\T)$ norm.
The families indexed by $j$ may be finite or countable.
We write $(h_j)\in L^q(\T;\ell^2)$ when
\[
 \left\|(h_j)\right\|_{L^q(\T;\ell^2)}
 :=\left\|\left(\sum_j|h_j|^2\right)^{1/2}\right\|_q<\infty.
\]
For an interval $I\subset\mathbb Z$, the operator $P_I$ is the bounded
Fourier projection on $L^q(\T)$ characterized by
\[
 \wh{P_Ih}(n)=\mathbf1_I(n)\wh h(n)\qquad(n\in\mathbb Z).
\]
The following estimates hold, with constants depending only on $q$.
\begin{enumerate}[label=(\roman*)]
\item Let $(I_j)$ be any family of intervals in $\mathbb Z$ and let
$(h_j)\in L^q(\T;\ell^2)$. Then
\[
 \left\|\left(\sum_j|P_{I_j}h_j|^2\right)^{1/2}\right\|_q
 \le C_q\left\|\left(\sum_j|h_j|^2\right)^{1/2}\right\|_q.
\]

\item Let $(I_j)$ be pairwise disjoint intervals in $\mathbb Z$.
If $2\le q<\infty$ and $h\in L^q(\T)$, then Rubio de Francia's
inequality gives
\[
 \left\|\left(\sum_j|P_{I_j}h|^2\right)^{1/2}\right\|_q
 \le C_q\|h\|_q.
\]
In the dual range $1<q\le2$, suppose that
$(u_j)\in L^q(\T;\ell^2)$ and that each $u_j$ has its Fourier
coefficients supported in $I_j$, that is,
\[
 \wh u_j(n)=0\qquad(n\notin I_j).
\]
Then
\[
 \left\|\sum_j u_j\right\|_q
 \le C_q\left\|\left(\sum_j|u_j|^2\right)^{1/2}\right\|_q.
\]
For a countable family, the series on the left converges
unconditionally in $L^q(\T)$.

\item For $h\in L^q(\T)$, define the one-sided Carleson maximal
operator by
\[
 (\mathcal C h)(e^{i\theta})
 :=\sup_{N\in\mathbb Z}
       |(P_{(-\infty,N]}h)(e^{i\theta})|.
\]
For a trigonometric polynomial $h$, this is
$\sup_{N\in\mathbb Z}|\sum_{n\le N}\wh h(n)e^{in\theta}|$;
the projection definition specifies its meaning for general
$h\in L^q(\T)$.
For every $(h_j)\in L^q(\T;\ell^2)$, one has
\[
 \left\|\left(\sum_j|\mathcal C h_j|^2\right)^{1/2}\right\|_q
 \le C_q\left\|\left(\sum_j|h_j|^2\right)^{1/2}\right\|_q.
\]
\end{enumerate}
\end{lemma}

\begin{proof}
We first work with finite families of trigonometric polynomials.

\smallskip\noindent\emph{(i) Interval projections.}
Let $P_+$ be the Riesz projection onto the nonnegative frequencies.
The scalar M.~Riesz theorem \cite[Chapter~III, \S2]{Garnett} gives its
$\ell^2$-valued bound by randomization. Indeed, for independent
Rademacher signs $(\epsilon_j)$, Khintchine's inequality
\cite[Appendix~C, \S C.2]{Grafakos} yields
\begin{align*}
 \left\|(P_+h_j)_j\right\|_{L^q(\T;\ell^2)}^q
 &\asymp_q
 \mathbb E\left\|P_+\Big(\sum_j\epsilon_jh_j\Big)\right\|_q^q\\
 &\lesssim_q\mathbb E\left\|\sum_j\epsilon_jh_j\right\|_q^q
 \asymp_q\left\|(h_j)_j\right\|_{L^q(\T;\ell^2)}^q.
\end{align*}
Write $(M_kh)(e^{i\theta})=e^{ik\theta}h(e^{i\theta})$.
For $I=[a,b]\cap\mathbb Z$ one has
\[
 P_I=M_aP_+M_{-a}-M_{b+1}P_+M_{-(b+1)},
\]
Half-line projections are treated in the same way, using either
$P_+$ or $\operatorname{Id}-P_+$.
Since coordinatewise modulations preserve the pointwise $\ell^2$ norm,
the preceding bound and the triangle inequality prove (i).

\smallskip\noindent\emph{(ii) Disjoint intervals and duality.}
The analysis estimate is the periodic form of Rubio de Francia's
theorem \cite[Theorem~1.2]{RdF}; the passage to the circle is recalled
below. To obtain the synthesis estimate, let $1<q\le2$ and put
$q'=q/(q-1)\ge2$. Since $P_{I_j}u_j=u_j$, for $g\in L^{q'}(\T)$
we have, with $dm$ denoting normalized arc measure,
\begin{align*}
 \left|\int_\T\Big(\sum_j u_j\Big)\overline g\,dm\right|
 &=\left|\sum_j\int_\T u_j\,\overline{P_{I_j}g}\,dm\right|\\
 &\le\left\|\Big(\sum_j|u_j|^2\Big)^{1/2}\right\|_q
       \left\|\Big(\sum_j|P_{I_j}g|^2\Big)^{1/2}\right\|_{q'}\\
 &\le C_{q'}\left\|\Big(\sum_j|u_j|^2\Big)^{1/2}\right\|_q
       \|g\|_{q'}.
\end{align*}
Taking the supremum over $\|g\|_{q'}\le1$ proves the claim.

\smallskip\noindent\emph{(iii) The maximal estimate.}
Choose $1<r<\min(q,2)$. For the real-line Carleson operator
$\mathcal C_{\mathbb R}$, \cite[Corollary~1.3 and \S2]{GMS} gives
\begin{align*}
 \left\|(\mathcal C_{\mathbb R}F_j)_j\right\|_{L^q(\mathbb R;\ell^2)}
 &\lesssim_{q,r}
 \left\|(M_rF_j)_j\right\|_{L^q(\mathbb R;\ell^2)}\\
 &\lesssim_{q,r}
 \left\|(F_j)_j\right\|_{L^q(\mathbb R;\ell^2)}.
\end{align*}
Here $M_rF=(M(|F|^r))^{1/r}$ and $M$ is the Hardy--Littlewood
maximal operator. The second inequality is Fefferman--Stein
\cite{FS}, applied with spatial exponent $q/r>1$ and sequence
exponent $2/r>1$. The one-sided Fourier maximal operator satisfies
the same bound, since half-line projections are linear combinations
of the identity and modulated Hilbert transforms.

To pass from the real line to the circle in (ii) and (iii), use the
localization argument underlying de Leeuw's restriction principle
\cite[Corollary~2.5]{deLeeuw}. Choose the real interval endpoints halfway between
integers and a nonzero Schwartz function $\psi$ with compactly
supported Fourier transform. For sufficiently small $\varepsilon>0$,
each projection of $\psi(\varepsilon t)h(e^{it})$ equals
$\psi(\varepsilon t)(P_Ih)(e^{it})$, because the multiplier is constant
on each frequency packet around an integer. Apply the real-line
estimate to these localized polynomials, multiply its $q$-th power
by $\varepsilon$, and let $\varepsilon\downarrow0$.
Periodic averaging cancels the common factor
$\int_{\mathbb R}|\psi|^q$, giving the circle estimate with the same
constant. In (iii), first take finitely many maximal cutoffs, placed
at half-integers, and then use monotone convergence.

Finally, density extends the bounds to finite families in
$L^q(\T;\ell^2)$; for the synthesis estimate one may use Fej\'er
means, which preserve the prescribed Fourier supports.
Monotone convergence gives the countable square-function bounds.
For countable synthesis, apply the finite-family estimate to tails:
their $L^q(\T;\ell^2)$ norms tend to zero, so the series converges
unconditionally in $L^q(\T)$.
\end{proof}

We will also use the nonnegative form of Rosenthal's inequality
\cite[equation~(1) and Theorem~1]{IS}: for independent nonnegative $X_j$ and $r>1$,
\begin{equation}\label{eq:rosenthal}
\mathbb E\left(\sum_jX_j\right)^r
\le C_r\left[\left(\sum_j\mathbb EX_j\right)^r+\sum_j\mathbb EX_j^r\right].
\end{equation}

\subsection{The square function and the Hankel representation}

\begin{theorem}[Square function]\label{thm:sq}
If $N_a$ is bounded on $H^2$, then for every $1<p<\infty$ and every polynomial $f$,
\begin{equation}\label{eq:sq}
 \|\cS f\|_{L^p(\T)}
 \le C_p\,\|N_a\|_{H^2\to H^2}\,\|f\|_{H^p}.
\end{equation}
\end{theorem}

\begin{proof}[Proof of Theorem~\ref{thm:sq}]
Let $\Gamma_\bh$ be the vector-valued Hankel operator with coefficients
$\bh_j=(\delta_m^{1/2}a_{m-j})_{m\ge0}$, where $a_l=0$ for $l<0$.
Lemma~\ref{lem:hankel} gives
\[
 \|\Gamma_\bh\|_{H^2\to H^2(\ell^2)}
 =\|N_a\|_{H^2\to H^2}.
\]
Moreover, $\sum_j\|\bh_j\|_{\ell^2}^2=\|N_a1\|_{H^2}^2<\infty$,
so $\bh(\zeta)=\sum_{j\ge0}\bh_j\zeta^j$ belongs to $H^2(\ell^2)$.

For a polynomial $f$, put $\check f(\zeta)=f(\bar\zeta)$ on $\T$.
Then $\Gamma_\bh f=P_+(\bh\,\check f)$, and its $m$-th coordinate is
\begin{align*}
 (\Gamma_\bh f)_m(\zeta)
 &=\delta_m^{1/2}\sum_{i=0}^m
       \wh{(\varphi f)}(m-i)\zeta^i\\[3pt]
 &=\delta_m^{1/2}\zeta^m
       [S_m(\varphi f)](\bar\zeta).
\end{align*}
Since $|\zeta^m|=1$, this identifies the square function pointwise:
\[
 \|(\Gamma_\bh f)(\zeta)\|_{\ell^2}^2
 =\sum_{m\ge0}\delta_m|S_m(\varphi f)(\bar\zeta)|^2
 =(\cS f(\bar\zeta))^2.
\]

\smallskip\noindent
To apply Page's theorem, put $\mathcal K=\mathbb C\oplus\ell^2$
and define $\widetilde\Gamma(f,F)=(0,\Gamma_\bh f)$ on
$H^2(\mathcal K)$. The dependence of the Hankel matrix entries on
$i+k$ gives $S_{\mathcal K}^*\widetilde\Gamma=\widetilde\Gamma S_{\mathcal K}$,
where $S_{\mathcal K}$ is multiplication by $z$. Page's representation
theorem \cite[Theorem~5, p.~534]{Page}
(see also \cite[Chapter~2, \S2]{Peller}) gives an essentially bounded
$\mathcal B(\mathcal K)$-valued symbol $\Psi$ such that
\[
 \widetilde\Gamma u=P_+(\Psi\,\check u),\qquad
 \|\Psi\|_{L^\infty}=\|\widetilde\Gamma\|=\|\Gamma_\bh\|,
\]
where $\check u(\zeta)=u(\bar\zeta)$. Let $\bg(\zeta)$ be the
$\ell^2$-coordinate of $\Psi(\zeta)(1,0)$. Applying the representation
to $(f,0)$ and taking its $\ell^2$-coordinate gives
\[
 \Gamma_\bh f=P_+(\bg\,\check f).
\]
We have $\|\bg\|_{L^\infty(\T;\ell^2)}\le\|\Psi\|_{L^\infty}$.
The reverse bound $\|\Gamma_\bh\|\le\|\bg\|_{L^\infty(\T;\ell^2)}$
follows from the last representation, since $P_+$ is contractive on
$L^2(\T;\ell^2)$ and reflection preserves the $L^2$ norm.
Taking $f=1$ also gives $P_+\bg=\Gamma_\bh 1=\bh$. Thus
\[
 \wh{\bg}(j)=\bh_j\quad(j\ge0),
 \qquad
 \|\bg\|_{L^\infty(\T;\ell^2)}=\|\Gamma_\bh\|.
\]
The coefficient identity prescribes only the nonnegative Fourier
coefficients of $\bg$. Page's theorem supplies the negative-frequency
part so that the full symbol is bounded. Since $\check f$ has only
nonpositive frequencies, $(\bg-\bh)\check f$ has strictly negative
frequencies and is removed by $P_+$. This allows us to use the bounded
symbol $\bg$ to estimate the same Hankel operator.

The Hilbert-valued Riesz bound from Lemma~\ref{lem:fourier}(i),
together with invariance of arc measure under $\zeta\mapsto\bar\zeta$,
now yields
\begin{align*}
 \|\cS f\|_{L^p(\T)}
 &=\|\Gamma_\bh f\|_{L^p(\T;\ell^2)}
   =\|P_+(\bg\check f)\|_{L^p(\T;\ell^2)}\\[3pt]
 &\le C_p\|\bg\check f\|_{L^p(\T;\ell^2)}
 \le C_p\|\bg\|_{L^\infty(\T;\ell^2)}\|f\|_{H^p}\\[3pt]
 &=C_p\|N_a\|_{H^2\to H^2}\|f\|_{H^p}.
\end{align*}
\end{proof}

\Needspace{14\baselineskip}
\begin{theorem}[Exponent independence of the square function]\label{thm:Sp}
Let $1<p<\infty$ and let $K_p\in[0,\infty]$ be the best constant such that
\[
 \|\cS f\|_{L^p}\le K_p\|f\|_{H^p}
 \qquad\text{for every analytic polynomial }f.
\]
Then
\[
\|N_a\|_{H^2\to H^2}\le K_p\le C_p\,\|N_a\|_{H^2\to H^2}.
\]
So $\cS$ is bounded from $H^p$ to $L^p$ for one $p$ iff for all $p$ iff $\tau(a)<\infty$.
\end{theorem}

\begin{proof}[Proof of Theorem~\ref{thm:Sp}]
The upper bound is Theorem~\ref{thm:sq}. For the lower bound, assume
$K=K_p<\infty$. Given $F\in H^p$, choose analytic polynomials
$f_j\to F$ in $H^p$. For each fixed $m$, convergence of the Taylor
coefficients gives $S_m(\varphi f_j)\to S_m(\varphi F)$ uniformly on
$\T$. Fatou's lemma, first for the sum over $m$ and then for the
integral on $\T$, yields
\[
 \|\cS F\|_{L^p}^p
 \le\liminf_{j\to\infty}\|\cS f_j\|_{L^p}^p
 \le K^p\|F\|_{H^p}^p.
\]
Thus the same square-function formula defines $\cS F\in L^p$ for
every $F\in H^p$, with the same bound $K$.

We use the same vector-valued Hankel operator $\Gamma_\bh$ as above,
now defined coordinatewise before any
$H^2$ boundedness is known. For $F\in H^p$, its coordinates are
\begin{align*}
 (\Gamma_\bh F)_m(\zeta)
 &=\delta_m^{1/2}\zeta^m[S_m(\varphi F)](\bar\zeta)\\[4pt]
 &=\delta_m^{1/2}\sum_{i=0}^{m}
       \wh{(\varphi F)}(m-i)\zeta^i.
\end{align*}
Thus $\|\Gamma_\bh F(\zeta)\|_{\ell^2}=\cS F(\bar\zeta)$.

For $G=(G_m)$ with finitely many nonzero coordinates $G_m\in H^{p'}$,
write $G_m^\sharp(z)=\overline{G_m(\bar z)}$. Comparing Fourier
coefficients gives
\begin{equation}\label{eq:bil}
 \begin{aligned}
  \sum_m\int_\T(\Gamma_\bh F)_m\,\overline{G_m}\,dm
  &=\sum_m\delta_m^{1/2}\sum_{s=0}^{m}
        a_{m-s}\,\wh{(FG_m^\sharp)}(s)\\[4pt]
  &=:\mathcal L(FG^\sharp).
 \end{aligned}
\end{equation}
The bilinear form therefore depends only on the product $FG^\sharp$
(the Hankel property), and
\begin{equation}\label{eq:bilbound}
 \begin{aligned}
  |\mathcal L(FG^\sharp)|
  &\le\|\cS F\|_p\,\|G\|_{L^{p'}(\ell^2)}\\[4pt]
  &\le K\|F\|_p\,\|G\|_{L^{p'}(\ell^2)}.
 \end{aligned}
\end{equation}

Now let $f$ and $G=(G_m)_{m=0}^{M}$ be analytic polynomials, with
$\Phi=fG^\sharp\ne0$, and set $w=\|\Phi\|_{\ell^2}$ on $\T$.
Then $\log w\in L^1$: a nonzero coordinate supplies an integrable
lower bound, while $w$ is bounded above. Choose the outer function
$F$ with $|F|=w^{1/p}$ on $\T$, and put $\Psi=\Phi/F$.
The outer-function construction and the scalar factorization facts
used here are in \cite[Chapter~II]{Garnett} and \cite[Chapter~2]{Duren}.

Each $\Psi_m$ lies in the Smirnov class and has boundary values in
$L^{p'}$, so Smirnov's theorem \cite[Theorem~2.11]{Duren} gives
$\Psi_m\in H^{p'}$. Moreover,
\[
 \|F\|_p^p
 =\|\Psi\|_{L^{p'}(\ell^2)}^{p'}
 =\|\Phi\|_{L^1(\ell^2)}.
\]
Apply \eqref{eq:bilbound} to $(F,\Psi^\sharp)$ to obtain
\begin{align*}
 |\mathcal L(\Phi)|
 &=|\mathcal L(F\Psi)|\\[4pt]
 &\le K\|\Phi\|_{L^1(\ell^2)}
 \le K\|f\|_2\,\|G\|_{L^2(\ell^2)}.
\end{align*}

By \eqref{eq:bil} with $F=f$, the left side is
$|\sum_{m=0}^{M}\langle(\Gamma_\bh f)_m,G_m\rangle|$.
Since these coordinates are analytic polynomials, taking the supremum
over $G$ gives
\[
 \sum_{m=0}^{M}\delta_m\|S_m(\varphi f)\|_2^2
 =\sum_{m=0}^{M}\|(\Gamma_\bh f)_m\|_2^2
 \le K^2\|f\|_2^2.
\]
Let $M\to\infty$ and apply Lemma~\ref{lem:hankel}:
\[
 \|N_af\|_2^2
 =\sum_{m\ge0}\delta_m\|S_m(\varphi f)\|_2^2
 \le K^2\|f\|_2^2.
\]
This proves $\|N_a\|_{H^2\to H^2}\le K$.
\end{proof}

\noindent The preceding argument is the scalar-outer-function proof
of the factorization $H^1(\ell^2)=H^p\cdot H^{p'}(\ell^2)$ needed here.
It uses scalar Smirnov theory coordinatewise, so its constants do not
depend on the number of coordinates.

\subsection{Sufficiency}

We now show that boundedness of $N_a$ on $H^2$ implies boundedness
on $H^p$. The block identity and the random cut points below are used
only in this subsection. Recall that
$A_n=\sum_{k=0}^na_k$ and that the doubling-time partition is given by
\[
 n_0=0,\qquad
 n_{j+1}=\min\{n:A_n\ge2A_{n_j}\},\qquad
 W_j=[n_j,n_{j+1})\cap\mathbb Z.
\]

\begin{lemma}[Blocks]\label{lem:block}
Let $f$ be a polynomial, and write
\[
 \varepsilon_m=\frac1{A_m}-\frac1{A_{m+1}}
 \quad(m\ge0).
\]
Set $S_{-1}(\varphi f)=0$. For integers $-1\le m'<m$, one has
\begin{equation}\label{eq:block}
 P_{(m',m]}N_af
 =\sum_{k=m'+1}^{m}\varepsilon_k S_k(\varphi f)
   +\frac{S_m(\varphi f)}{A_{m+1}}
   -\frac{S_{m'}(\varphi f)}{A_{m'+1}}.
\end{equation}
Here $P_{(m',m]}$ retains precisely the coefficients of degrees
$m'+1,\ldots,m$.

If $m'\in W_{j-1}$ and $m\in W_j$ for some $j\ge1$, then, pointwise
on $\T$,
\begin{equation}\label{eq:blockbound}
 |P_{(m',m]}N_af|^2
 \le12\sum_{k=m'+1}^{m}\delta_k|S_k(\varphi f)|^2
    +3\frac{|S_m(\varphi f)|^2}{A_m^2}
    +3\frac{|S_{m'}(\varphi f)|^2}{A_{m'}^2}.
\end{equation}
For the initial block, $m'=-1$ and $m\in W_0$, the corresponding bound is
\[
 |S_m(N_af)|^2
 \le12\sum_{k=0}^{m}\delta_k|S_k(\varphi f)|^2
    +3\frac{|S_m(\varphi f)|^2}{A_m^2}.
\]

Finally, for every $m\ge0$,
\[
 \frac{|S_m(\varphi f)|}{A_m}
 \le2\max_{0\le k\le m}|S_k(N_af)|
 \qquad\text{on }\T.
\]
\end{lemma}

\begin{proof}
For $m'<n\le m$, telescoping gives
\[
 \frac1{A_n}=\frac1{A_{m+1}}+\sum_{k=n}^{m}\varepsilon_k.
\]
Multiply by $\wh{(\varphi f)}(n)z^n$ and sum over $n$. Interchanging
the finite sums yields
\[
 P_{(m',m]}N_af
 =\frac{S_m(\varphi f)-S_{m'}(\varphi f)}{A_{m+1}}
   +\sum_{k=m'+1}^{m}\varepsilon_k
      \bigl(S_k(\varphi f)-S_{m'}(\varphi f)\bigr),
\]
which simplifies to \eqref{eq:block}.

Since $\varepsilon_k/A_k\le\delta_k$, Cauchy--Schwarz gives
\[
 \left|\sum_{k=m'+1}^{m}\varepsilon_k S_k(\varphi f)\right|^2
 \le\left(\sum_{k=m'+1}^{m}\varepsilon_kA_k\right)
       \sum_{k=m'+1}^{m}\delta_k|S_k(\varphi f)|^2.
\]
Now $\varepsilon_kA_k=1-A_k/A_{k+1}\le\log(A_{k+1}/A_k)$.
The terms with $k,k+1$ in one block sum to at most $\log2$.
For two adjacent blocks there are at most two further terms, each at
most $1$, so
\[
 \sum_{k=m'+1}^{m}\varepsilon_kA_k\le2\log2+2\le4.
\]
Apply $|u+v+w|^2\le3(|u|^2+|v|^2+|w|^2)$ to \eqref{eq:block},
using $A_{m+1}\ge A_m$ and $A_{m'+1}\ge A_{m'}$, to obtain
\eqref{eq:blockbound}. For $m'=-1$, the last term in \eqref{eq:block}
vanishes, and the same argument gives the initial-block bound.

Finally, summation by parts gives the explicit identity
\[
 \frac{S_m(\varphi f)}{A_m}
 =S_m(N_af)
   -\sum_{k=0}^{m-1}\frac{A_{k+1}-A_k}{A_m}\,S_k(N_af).
\]
The nonnegative coefficients in the sum have total
$(A_m-A_0)/A_m\le1$, proving the last assertion.
\end{proof}

\medskip
\noindent\textbf{Random cut points.}
Fix a polynomial $f$.
Choose an integer $J\ge1$ and put $M=n_{J+1}-1$.
For $0\le j<J$, define the positive mass
\[
 \delta(W_j):=\sum_{m\in W_j}\delta_m
 =A_{n_j}^{-2}-A_{n_{j+1}}^{-2}
 \ge\frac34A_{n_j}^{-2}.
\]
Choose independent random integers $m_j\in W_j$ with distribution
\[
 \mathbb P(m_j=m)=\frac{\delta_m}{\delta(W_j)}
 \qquad(m\in W_j,\ 0\le j<J).
\]
Set $m_{-1}=-1$ and $m_J=M$, and put
\[
 I_j=(m_{j-1},m_j]\cap\mathbb Z
 \qquad(0\le j\le J).
\]
These intervals partition $[0,M]\cap\mathbb Z$, so
$S_M(N_af)=\sum_{j=0}^{J}P_{I_j}N_af$.

For $0\le j<J$, let $\beta_j=|S_{m_j}(\varphi f)|/A_{m_j}$, viewed
pointwise on $\T$. Lemma~\ref{lem:block} and the above distribution give
\begin{equation}\label{eq:sumu}
 \begin{aligned}
  \sum_{j=0}^{J}|P_{I_j}N_af|^2
  &\le12(\cS f)^2+6\sum_{j=0}^{J-1}\beta_j^2
       +3\frac{|S_M(\varphi f)|^2}{A_M^2},\\[6pt]
  \E\,\beta_j^2
  &=\frac1{\delta(W_j)}
       \sum_{m\in W_j}\delta_m\frac{|S_m(\varphi f)|^2}{A_m^2}
   \le\frac43\sum_{m\in W_j}\delta_m|S_m(\varphi f)|^2.
 \end{aligned}
\end{equation}
The second bound uses $A_m\ge A_{n_j}$ for $m\in W_j$.
The deterministic endpoint is controlled, for $1<p<\infty$, by
\[
 \|S_M(\varphi f)\|_{L^p(\T)}
 =\|S_M((S_M\varphi)f)\|_{L^p(\T)}
 \le B_pA_M\|f\|_{H^p},
\]
where $B_p=\sup_{\nu\ge0}\|S_\nu\|_{H^p\to H^p}$; here
$\|S_M\varphi\|_\infty\le A_M$.

\medskip
The following two estimates explain the roles of the two sides of $p=2$.
For $p\le2$, no condition on the frequency partition is needed. For
$p>2$, Theorem~\ref{thm:LP} will supply the stated hypothesis whenever
$\tau(a)<\infty$.

\Needspace{8\baselineskip}
\begin{theorem}[$H^p$, $1<p\le2$, no doubling]\label{thm:low}
If $N_a$ is bounded on $H^2$ (equivalently $\tau(a)<\infty$), then $N_a$ is bounded on $H^p$ for
every $1<p\le2$, and $\|N_a\|_{H^p\to H^p}\le C_p\|N_a\|_{H^2\to H^2}$.
\end{theorem}

\begin{proof}[Proof of Theorem~\ref{thm:low}]
Fix $1<p\le2$ and a polynomial $f$, and use the random cut points
introduced above. The projections $P_{I_j}N_af$ have disjoint
Fourier supports, so Lemma~\ref{lem:fourier}(ii) gives
\[
 \|S_M(N_af)\|_p
 \le C_p\left\|\left(\sum_{j=0}^{J}|P_{I_j}N_af|^2\right)^{1/2}\right\|_p.
\]

Since $t\mapsto t^{p/2}$ is concave, Jensen's inequality and
\eqref{eq:sumu} yield
\begin{align*}
 \E\left\|\left(\sum_{j=0}^{J-1}\beta_j^2\right)^{1/2}\right\|_p^p
 &\le\int_\T
       \left(\sum_{j=0}^{J-1}\E\,\beta_j^2\right)^{p/2}\,dm\\[4pt]
 &\le C_p\|\cS f\|_p^p.
\end{align*}
Thus some choice of the cuts satisfies
$\|(\sum_{j=0}^{J-1}\beta_j^2)^{1/2}\|_p\le C_p\|\cS f\|_p$.
For this choice, \eqref{eq:sumu} and the deterministic endpoint bound give
\begin{align*}
 \|S_M(N_af)\|_p
 &\le C_p\bigl(\|\cS f\|_p+\|f\|_{H^p}\bigr)\\[4pt]
 &\le C_p\|N_a\|_{H^2\to H^2}\|f\|_{H^p},
\end{align*}
where the last step uses Theorem~\ref{thm:sq} and
$\|N_a\|_{H^2\to H^2}\ge1$.

The bound is independent of $J$. As $J\to\infty$, reflexivity of $H^p$
gives a weakly convergent subsequence of these Taylor sections.
Their coefficients identify the limit as $N_af$.
Density of polynomials completes the proof.
\end{proof}

\Needspace{8\baselineskip}
\begin{theorem}[$H^p$, $p>2$]\label{thm:high}
Let $2<p<\infty$. If $N_a$ is bounded on $H^2$ and the partition $\{W_j\}$ satisfies
$(\mathrm{LP}_{p'})$, then $N_a$ is bounded on $H^p$, with norm $\le C(p,L)\,\|N_a\|_{H^2\to H^2}$.
\end{theorem}

\begin{proof}[Proof of Theorem~\ref{thm:high}]
Fix a polynomial $f$, let $p>2$ and $q=p'$, and set
\[
 X=\|S_M(N_af)\|_p,\qquad
 Y=\|\cS f\|_p+\|f\|_{H^p}.
\]
All constants below are independent of $J$ and of the cut points.

Fix a choice of cuts and a trigonometric polynomial $g$ with $\|g\|_q=1$.
For $j\ge1$, $I_j\subset W_{j-1}\cup W_j$, while
$I_0\subset W_0$. Hence
\[
 P_{I_0}g=P_{I_0}P_{W_0}g,\qquad
 P_{I_j}g=P_{I_j}(P_{W_{j-1}}g+P_{W_j}g)\quad(1\le j\le J).
\]
Lemma~\ref{lem:fourier}(i) and $(\mathrm{LP}_q)$ give
\begin{align*}
 \left\|\left(\sum_{j=0}^{J}|P_{I_j}g|^2\right)^{1/2}\right\|_q
 &\le C_q\left\|\left(
 |P_{W_0}g|^2+\sum_{j=1}^{J}|P_{W_{j-1}}g+P_{W_j}g|^2
 \right)^{1/2}\right\|_q\\[3pt]
 &\le 2C_q\left\|\left(\sum_{j=0}^{J}|P_{W_j}g|^2\right)^{1/2}\right\|_q\\[3pt]
 &\le 2C_qL\|g\|_q=2C_qL.
\end{align*}

\Needspace{10\baselineskip}
The Fourier support of $P_{I_j}N_af$ lies in $I_j$, so its pairing
with $g$ equals its pairing with $P_{I_j}g$. Using Cauchy--Schwarz in $j$,
H\"older's inequality with $1/p+1/q=1$, and the preceding estimate, we obtain
\begin{align*}
 \big|\langle S_M(N_af),g\rangle\big|
 &=\left|\sum_{j=0}^{J}\langle P_{I_j}N_af,P_{I_j}g\rangle\right|\\[3pt]
 &\le\int_\T\sum_{j=0}^{J}|P_{I_j}N_af|\,|P_{I_j}g|\,dm\\[3pt]
 &\le\int_\T\left(\sum_{j=0}^{J}|P_{I_j}N_af|^2\right)^{1/2}
            \left(\sum_{j=0}^{J}|P_{I_j}g|^2\right)^{1/2}\,dm\\[3pt]
 &\le\left\|\left(\sum_{j=0}^{J}|P_{I_j}N_af|^2\right)^{1/2}\right\|_p
      \left\|\left(\sum_{j=0}^{J}|P_{I_j}g|^2\right)^{1/2}\right\|_q\\[3pt]
 &\le 2C_qL\left\|\left(\sum_{j=0}^{J}|P_{I_j}N_af|^2\right)^{1/2}\right\|_p.
\end{align*}
The left-hand side is independent of the cuts. Averaging over them and
taking the supremum over $g$ therefore gives, by $L^p$ duality,
\[
 X\le 2C_qL\,\E
       \left\|\left(\sum_{j=0}^{J}|P_{I_j}N_af|^2\right)^{1/2}\right\|_p.
\]

To estimate this expectation, apply Rosenthal's inequality
\eqref{eq:rosenthal} pointwise with $r=p/2$:
\[
 \E\left(\sum_{j=0}^{J-1}\beta_j^2\right)^r
 \le C_r\left[
       \left(\sum_{j=0}^{J-1}\E\,\beta_j^2\right)^r
       +\sum_{j=0}^{J-1}\E\,\beta_j^{2r}\right].
\]
Lemma~\ref{lem:block} gives $\beta_j\le2\mathcal C(S_M(N_af))$ for every
cut. The second-moment bound in \eqref{eq:sumu} therefore yields
\begin{align*}
 \sum_{j=0}^{J-1}\E\,\beta_j^p
 &\le \bigl[2\mathcal C(S_M(N_af))\bigr]^{p-2}
       \sum_{j=0}^{J-1}\E\,\beta_j^2\\[4pt]
 &\le C_p\bigl[\mathcal C(S_M(N_af))\bigr]^{p-2}(\cS f)^2.
\end{align*}
H\"older's inequality and the scalar Carleson--Hunt estimate
(Lemma~\ref{lem:fourier}(iii) and \cite[\S2]{GMS}) yield
\begin{align*}
 \int_\T\sum_{j=0}^{J-1}\E\,\beta_j^p\,dm
 &\le C_p\|\mathcal C(S_M(N_af))\|_p^{p-2}\|\cS f\|_p^2\\[4pt]
 &\le C_pX^{p-2}\|\cS f\|_p^2.
\end{align*}

Combining these bounds with \eqref{eq:sumu} and the deterministic
endpoint estimate, and using Jensen's inequality for the expectation,
we obtain
\[
 X\le C(p,L)\bigl(Y+X^{1-2/p}Y^{2/p}\bigr).
\]
Young's inequality absorbs the mixed term:
\[
 C(p,L)X^{1-2/p}Y^{2/p}
 \le\frac X2+C'(p,L)Y.
\]
Consequently, Theorem~\ref{thm:sq} and
$\|N_a\|_{H^2\to H^2}\ge1$ give
\begin{align*}
 \|S_M(N_af)\|_p
 &\le C(p,L)\bigl(\|\cS f\|_p+\|f\|_{H^p}\bigr)\\[4pt]
 &\le C(p,L)\|N_a\|_{H^2\to H^2}\|f\|_{H^p}.
\end{align*}

This bound is uniform in $J$. Reflexivity gives a weakly convergent
subsequence of the Taylor sections, whose coefficients identify the
limit as $N_af$. Extend by density of polynomials in $H^p$.
\end{proof}

\subsection{Necessity}\label{sec:nec}

We now prove necessity: boundedness of $N_a$ on $H^p$ implies boundedness
on $H^2$. By Theorem~\ref{thm:Sp}, it suffices to prove
\[
 \|\cS f\|_{L^p}\le C\,\|N_af\|_{H^p}
\]
for polynomials $f$. When $1<p<2$, the estimate below uses
$(\mathrm{LP}_p)$; this hypothesis is supplied in the completion
of the proof.

To express $\cS f$ in terms of $N_af$, define one further square function:
\begin{equation}\label{eq:Q}
 \begin{aligned}
  \cQ h&:=\left(\sum_{m\ge0}\delta_m
       \left|\sum_{n=0}^{m}A_n\wh h(n)z^n\right|^2\right)^{1/2},
       \\[6pt]
  \cS f&=\cQ(N_af).
 \end{aligned}
\end{equation}
The identity follows from
$\wh{(\varphi f)}(n)=A_n\wh{(N_af)}(n)$.
Thus $\cQ$ depends only on $(A_n)$. In particular,
$\|\cQ h\|_2=\|h\|_2$, since
$\sum_{m\ge n}\delta_m=A_n^{-2}$.

\Needspace{10\baselineskip}
\begin{theorem}[The range square function]\label{thm:Q}
Let $(A_n)$ be positive and nondecreasing, with $A_n\to\infty$ and doubling-time partition $\{W_j\}$.
\begin{enumerate}
\item[(i)] If $2\le p<\infty$, then $\|\cQ h\|_{L^p}\le C_p\|h\|_{H^p}$ for all $h\in H^p$.
\item[(ii)] If $1<p<2$ and $\{W_j\}$ satisfies $(\mathrm{LP}_p)$ with constant $L$, then
$\|\cQ h\|_{L^p}\le C_pL\|h\|_{H^p}$.
\end{enumerate}
\end{theorem}

\begin{proof}[Proof of Theorem~\ref{thm:Q}]
We first work with polynomials and finite truncations of the block sums.
All estimates will be uniform in the truncation.
Within this proof, normalize the weighted Taylor sections by putting
\[
 U_mh:=\frac1{A_m}\sum_{n=0}^{m}A_n\wh h(n)z^n,\qquad
 \omega_m:=\delta_mA_m^2=1-\frac{A_m^2}{A_{m+1}^2}.
\]
Then $(\cQ h)^2=\sum_{m\ge0}\omega_m|U_mh|^2$.

For an interval $[\alpha,\beta]$, summation of the increments of $A_n$
gives
\begin{equation}\label{eq:convex}
 \sum_{n=\alpha}^{\beta}A_n\wh h(n)z^n
 =A_\alpha P_{[\alpha,\beta]}h
  +\sum_{k=\alpha+1}^{\beta}(A_k-A_{k-1})P_{[k,\beta]}h.
\end{equation}
After division by $A_\beta$, the right side is a convex combination
of interval projections. For the current block and its endpoint, write
\[
 V_{j,m}:=\frac1{A_m}\sum_{n=n_j}^{m}A_n\wh h(n)z^n
 \quad(m\in W_j),\qquad
 Y_j:=V_{j,n_{j+1}-1}.
\]
Jensen's inequality therefore yields
\begin{equation}\label{eq:jensen}
 \begin{aligned}
  |V_{j,m}|^2
  &\le\frac{A_{n_j}}{A_m}|P_{[n_j,m]}h|^2\\[4pt]
  &\quad+\sum_{k=n_j+1}^{m}
       \frac{A_k-A_{k-1}}{A_m}|P_{[k,m]}h|^2.
 \end{aligned}
\end{equation}
The coefficients sum to one. The same bound applies to $Y_j$ by
taking $m=n_{j+1}-1$.

For $m\in W_j$, decompose
\[
 U_mh=V_{j,m}+\sum_{i<j}\frac{A_{n_{i+1}-1}}{A_m}Y_i.
\]
The doubling-time construction gives
\[
 \frac{A_{n_{i+1}-1}}{A_m}
 <\frac{2A_{n_i}}{A_{n_j}}\le2^{1-(j-i)},
 \qquad
 \sum_{m\in W_j}\omega_m
 \le4A_{n_j}^2\sum_{m\in W_j}\delta_m\le4.
\]
Cauchy--Schwarz gives
\[
 \left(\sum_{i<j}2^{1-(j-i)}|Y_i|\right)^2
 \le2\sum_{i<j}2^{1-(j-i)}|Y_i|^2.
\]
Summing over $m$ and then over $j$ yields the common pointwise bound
\begin{equation}\label{eq:Qpt}
 (\cQ h)^2
 \le2\sum_j\sum_{m\in W_j}\omega_m|V_{j,m}|^2
    +32\sum_i|Y_i|^2.
\end{equation}

For a fixed block $W_j$, the two sums in \eqref{eq:Qpt} contribute
$2\sum_{m\in W_j}\omega_m|V_{j,m}|^2+32|Y_j|^2$.
To represent these weights by probabilities, choose a right endpoint
$m\in W_j$ with probability $2\omega_m/40$, adding $32/40$ to the
probability of $m=n_{j+1}-1$ to account for $Y_j$.
The total probability of choosing an endpoint is at most
$(2\cdot4+32)/40=1$. With the remaining probability, take
$K_j=\varnothing$ and set $P_{\varnothing}=0$.

Given a chosen endpoint $m$, take $K_j=[n_j,m]$ with probability
$A_{n_j}/A_m$, or take $K_j=[k,m]$ with probability
$(A_k-A_{k-1})/A_m$ for $n_j<k\le m$. By \eqref{eq:jensen}, this
conditional average of $|P_{K_j}h|^2$ is at least $|V_{j,m}|^2$.
Since $Y_j=V_{j,n_{j+1}-1}$, averaging over the right endpoint yields
\[
 2\sum_{m\in W_j}\omega_m|V_{j,m}|^2+32|Y_j|^2
 \le40\,\E|P_{K_j}h|^2.
\]
Summing over the blocks in \eqref{eq:Qpt} gives
\[
 (\cQ h)^2\le40\,\E\sum_j|P_{K_j}h|^2.
\]

\smallskip\noindent\emph{(i) The range $p\ge2$.}
For every selection, the intervals are pairwise disjoint.
Lemma~\ref{lem:fourier}(ii) and Minkowski's inequality in $L^{p/2}$
therefore give
\[
 \|\cQ h\|_p^2
 \le40\,\E\left\|\sum_j|P_{K_j}h|^2\right\|_{p/2}
 \le40 C_p^2\|h\|_p^2.
\]

\smallskip\noindent\emph{(ii) The range $1<p<2$.}
For every interval $I\subset W_j$,
\[
 |P_Ih|\le2\,\mathcal C(P_{W_j}h),
\]
where $\mathcal C$ is the maximal operator from
Lemma~\ref{lem:fourier}(iii). Applying this bound to each $K_j$ in the
average gives
\[
 (\cQ h)^2\le160\sum_j|\mathcal C(P_{W_j}h)|^2.
\]
The vector-valued maximal estimate and $(\mathrm{LP}_p)$ give
\begin{align*}
 \|\cQ h\|_p
 &\le C_p\left\|\left(\sum_j|P_{W_j}h|^2\right)^{1/2}\right\|_p\\[4pt]
 &\le C_pL\|h\|_p.
\end{align*}
Finally, Fatou's lemma removes the truncation, and density extends
both estimates to $h\in H^p$.
\end{proof}

\Needspace{10\baselineskip}
\begin{theorem}[Necessity]\label{thm:nec}
Let $a_n\ge0$.
\begin{enumerate}
\item[(i)] If $N_a$ is bounded on $H^p$ for some $2\le p<\infty$, then $\tau(a)<\infty$, and
$\|N_a\|_{H^2\to H^2}\le C_p\|N_a\|_{H^p\to H^p}$.
\item[(ii)] If $N_a$ is bounded on $H^p$ for some $1<p<2$ and $\{W_j\}$ satisfies $(\mathrm{LP}_p)$,
then $\tau(a)<\infty$, and $\|N_a\|_{H^2\to H^2}\le C_pL\,\|N_a\|_{H^p\to H^p}$.
\end{enumerate}
\end{theorem}

\begin{proof}[Proof of Theorem~\ref{thm:nec}]
For a polynomial $f$, identity \eqref{eq:Q} and
Theorem~\ref{thm:Q} give
\begin{align*}
 \|\cS f\|_p
 &=\|\cQ(N_af)\|_p\\[4pt]
 &\le C\|N_af\|_{H^p}\\[4pt]
 &\le C\|N_a\|_{H^p\to H^p}\|f\|_{H^p},
\end{align*}
where $C=C_p$ in (i) and $C=C_pL$ in (ii).
Theorem~\ref{thm:Sp} therefore yields
\[
 \|N_a\|_{H^2\to H^2}
 \le K_p
 \le C\|N_a\|_{H^p\to H^p},
\]
and in particular $\tau(a)<\infty$.
\end{proof}

\subsection{Completion of the proof}

\begin{proof}[Proof of Theorem~\ref{thm:main}]
First assume $A_n\to\infty$. The equivalence (ii)$\Leftrightarrow$(iii)
is Theorem~\ref{thm:l2}.
If (iii) holds, Theorem~\ref{thm:low} gives boundedness for $p\le2$.
For $p>2$, Proposition~\ref{prop:nec} gives $\kappa_0(a)\le4\tau(a)$;
Theorem~\ref{thm:LP} supplies $(\mathrm{LP}_{p'})$, and
Theorem~\ref{thm:high} gives boundedness for $p>2$.
Conversely, boundedness on $H^p$ implies (iii) for $p\ge2$ by Theorem~\ref{thm:nec}(i).
For $p<2$, Proposition~\ref{prop:nec} bounds $\kappa_0(a)$ in terms of $p$
and $\|N_a\|_{H^p}$, so Theorem~\ref{thm:LP} supplies
$(\mathrm{LP}_p)$ and Theorem~\ref{thm:nec}(ii) applies.
All constants in this argument have the dependence stated in the theorem.

\smallskip\noindent\emph{Finite total mass and regularization.}
Now suppose $A_\infty<\infty$. Then $\varphi\in H^\infty$ with
$\|\varphi\|_\infty=A_\infty$, and Abel summation for the decreasing
sequence $(1/A_n)$ gives
\[
 \|N_a\|_{H^p\to H^p}
 \le B_p\,\frac{A_\infty}{a_0}
 \qquad(1<p<\infty).
\]
Thus boundedness itself is automatic. For the quantitative bounds,
fix $\epsilon>0$ and put
\[
 a_n^{(\epsilon)}=a_n+\epsilon,\qquad
 A_n^{(\epsilon)}=A_n+\epsilon(n+1)\longrightarrow\infty.
\]
The infinite-sum case therefore applies to every regularized sequence.
With $C$ denoting the ordinary Ces\`aro operator, the coefficientwise
identity is
\[
 \begin{aligned}
  N_{a^{(\epsilon)}}
  &=\Lambda_{m^{(\epsilon)}}N_a
    +\Lambda_{1-m^{(\epsilon)}}C,\\[5pt]
  m_n^{(\epsilon)}
  &=\frac{A_n}{A_n+\epsilon(n+1)}.
 \end{aligned}
\]
These two multipliers are bounded on $H^p$ uniformly in $\epsilon$
and $a$. Indeed, $0\le m_n^{(\epsilon)}\le1$ and
\[
 (m_n^{(\epsilon)}-m_{n+1}^{(\epsilon)})_+
 \le\frac{A_n\epsilon}
 {(A_n+\epsilon(n+1))(A_n+\epsilon(n+2))}
 \le\frac1{n+1}.
\]
On every dyadic interval, the sum of the downward increments is
uniformly bounded. The sum of the upward increments differs from it
by at most one, so the total variation is uniformly bounded as well.
To apply the real-line Marcinkiewicz theorem, extend
$m^{(\epsilon)}$ to $\mathbb R$ by smooth monotone interpolation on
each $[n,n+1]$, with zero derivatives at the endpoints, and extend
it constantly to the left of zero. Its variation on each positive
dyadic interval with integer endpoints is exactly the sum of the
corresponding discrete increments in absolute value. On the other
dyadic intervals its variation is at most one. The same bounds hold
for the complementary symbol. The real-line Marcinkiewicz theorem
\cite[Theorem~6.2.2]{Grafakos} and de Leeuw's restriction theorem
\cite[Corollary~2.5]{deLeeuw} give the required uniform $H^p$ bounds.
Using also the boundedness of $C$ on $H^p$, we obtain
\begin{equation}\label{eq:regularized-Hp}
 \|N_{a^{(\epsilon)}}\|_{H^p\to H^p}
 \le C_p\bigl(\|N_a\|_{H^p\to H^p}+1\bigr).
\end{equation}

Apply the coefficientwise identity to $D_J$. For $n\ge J$, the
Ces\`aro row sum is $J/(n+1)$, so \eqref{eq:rowsum} gives
\[
 R_J^{(\epsilon)}(n)
 =m_n^{(\epsilon)}R_J(n)+(1-m_n^{(\epsilon)})\frac{J}{n+1}.
\]
By convexity, $R_J^{(\epsilon)}(n)^2\le R_J(n)^2+J^2/(n+1)^2$.
Since $J\sum_{n\ge J}(n+1)^{-2}\le1$, we get
\[
 \tau(a^{(\epsilon)})\le\tau(a)+1.
\]

The infinite-mass result and \eqref{eq:regularized-Hp} bound
$\tau(a^{(\epsilon)})$ in terms of $p$ and $\|N_a\|_{H^p\to H^p}$.
When $1<p<2$, Proposition~\ref{prop:nec} first bounds
$\kappa_0(a^{(\epsilon)})$ uniformly using \eqref{eq:regularized-Hp},
and Theorem~\ref{thm:LP} supplies a uniform $(\mathrm{LP}_p)$ constant.
For sufficiency, $\tau(a^{(\epsilon)})\le\tau(a)+1$ bounds
$\|N_{a^{(\epsilon)}}\|_{H^p\to H^p}$ in terms of $p$ and $\tau(a)$.

For each $J$ and $n$, $R_J^{(\epsilon)}(n)\to R_J(n)$ as
$\epsilon\downarrow0$. Fatou's lemma and the supremum over $J$ give
$\tau(a)\le\liminf_{\epsilon\downarrow0}\tau(a^{(\epsilon)})$.
For each polynomial $f$, the Taylor coefficients of
$N_{a^{(\epsilon)}}f$ converge to those of $N_af$ and are bounded by
$\max_k|\wh f(k)|$ uniformly in $\epsilon$. Thus the functions converge
uniformly on compact subsets of $\D$. Taking norms on each circle of
radius $r<1$, then the supremum over $r$, and using polynomial density
gives
\[
 \|N_a\|_{H^p\to H^p}
 \le\liminf_{\epsilon\downarrow0}\|N_{a^{(\epsilon)}}\|_{H^p\to H^p}.
\]
These two limit estimates give the asserted bounds for $a$.
\end{proof}

\section{The scalar symbol: a general criterion and BMOA}\label{sec:bmo}

We express the testing condition entirely in terms of the symbol $G_a$
from \eqref{eq:symbol}. The truncation in \eqref{eq:truncated-functional}
retains the saturation of the row sums $R_J(n)$, and therefore requires
no doubling assumption.

\subsection{The general symbol criterion}

For integers $J,k\ge1$, write
\[
 b_k(J)=\sum_{j=kJ+1}^{(k+1)J}c_j
\]
for the sum of the coefficients of $G_a$ over the $k$th block of length $J$
in \eqref{eq:truncated-functional}.

\begin{theorem}[Truncated coefficient criterion]\label{thm:symbol}
Let $a_n\ge0$, $a_0>0$, and let $\mathcal F_{\mathrm{tr}}(G_a)$ be
defined by \eqref{eq:truncated-functional}. Then
\begin{equation}\label{eq:truncated-comparison}
 \tau(a)\asymp 1+\mathcal F_{\mathrm{tr}}(G_a),
\end{equation}
with absolute comparison constants. Consequently, for every $1<p<\infty$,
\begin{equation}\label{eq:symbol-characterization}
 N_a:H^p\longrightarrow H^p\text{ is bounded}
 \quad\Longleftrightarrow\quad
 \mathcal F_{\mathrm{tr}}(G_a)<\infty.
\end{equation}
Moreover,
\[
 \|N_a\|_{H^2\to H^2}^{\,2}\asymp
 1+\mathcal F_{\mathrm{tr}}(G_a).
\]
\end{theorem}

\begin{proof}
For $n\ge J\ge1$, put
\[
 s_J(n)=\sum_{j=n-J+1}^nc_j.
\]
Since $1-c_j=A_{j-1}/A_j$ for $j\ge1$, we have the exact identity
\[
 R_J(n)=1-\frac{A_{n-J}}{A_n}
       =1-\prod_{j=n-J+1}^n(1-c_j).
\]
The inequalities $1-t\le e^{-t}$ and
$1-\prod_j(1-t_j)\le\sum_jt_j$, for $0\le t_j\le1$, imply
\begin{equation}\label{eq:truncated-window}
 (1-e^{-1})\min\{1,s_J(n)\}\le R_J(n)\le\min\{1,s_J(n)\}.
\end{equation}
Thus, if
\[
 T=\sup_{J\ge1}\frac1J\sum_{n\ge J}\min\{1,s_J(n)\}^2,
\]
then
\[
 1+(1-e^{-1})^2T\le\tau(a)\le1+T.
\]

It remains to compare sliding windows with disjoint blocks.
For $(k+1)J<n\le(k+2)J$ we have
$s_J(n)\le b_k(J)+b_{k+1}(J)$. Applying
\[
 \min\{1,x+y\}^2\le2\min\{1,x\}^2+2\min\{1,y\}^2
 \qquad(x,y\ge0)
\]
and summing over these windows gives
\[
 \frac1J\sum_{n>2J}\min\{1,s_J(n)\}^2
 \le4\sum_{k\ge1}\min\{1,b_k(J)\}^2.
\]
The terms $J\le n\le2J$ contribute at most $(J+1)/J\le2$.
Consequently,
\[
 T\le2+4\mathcal F_{\mathrm{tr}}(G_a).
\]

Conversely, for every integer $(k+1)J<n\le(k+2)J$,
the window $[n-2J+1,n]$ contains the block $[kJ+1,(k+1)J]$.
Therefore
\[
 J\sum_{k\ge1}\min\{1,b_k(J)\}^2
 \le\sum_{n>2J}\min\{1,s_{2J}(n)\}^2
 \le2JT.
\]
Taking the supremum over $J$ yields
$\mathcal F_{\mathrm{tr}}(G_a)\le2T$ and proves
\eqref{eq:truncated-comparison}. The remaining assertions follow from
Theorem~\ref{thm:main}.
\end{proof}

\subsection{Doubling partial sums and BMOA}

Recall that $\BMOA$ consists of the analytic functions whose boundary
values have bounded mean oscillation. The Poisson-oscillation
characterization \cite[Chapter~VI, \S\S1--3]{Garnett} gives the equivalent norm
\[
 \|g\|_{\BMOA}=|g(0)|+
 \sup_{w\in\D}\left(\int_{\T}|g(\zeta)-g(w)|^2
                  P_w(\zeta)\,dm(\zeta)\right)^{1/2}.
\]
For $g(z)=\sum_{n\ge0}d_nz^n$ with $d_n\ge0$, Fefferman's
coefficient criterion states that
\begin{equation}\label{eq:fefferman}
 g\in\BMOA\quad\Longleftrightarrow\quad
 \mathcal F(g):=\sup_{J\ge1}\sum_{k\ge1}
       \left(\sum_{j=kJ+1}^{(k+1)J}d_j\right)^2<\infty.
\end{equation}
For the criterion, see Garc\'ia-Cuerva, Kazarian, Kolyada and Torrea
\cite[Theorem~9.1, equation~(9.8), and Corollary~9.2 with the paragraph
following it]{GCKKT}. For an early account of the corresponding $H^1$
multiplier theorem, see Anderson and Shields \cite{AnSh}.
The source uses the blocks $[kJ,(k+1)J-1]$. Each block in
\eqref{eq:fefferman} lies in the union of two consecutive source
blocks. Conversely, each source block with $k\ge2$ lies in the union
of two consecutive blocks in \eqref{eq:fefferman}. For $k=1$ the
only extra coefficient is $d_J$, which is bounded uniformly in $J$
by the $J=1$ test together with $d_1$. By positivity and
$(u+v)^2\le2(u^2+v^2)$, the two conditions are equivalent.
We also use the classical Hankel characterization:
$\Gamma_g=(d_{j+k})_{j,k\ge0}$ is bounded on $\ell^2$ if and only if
$g\in\BMOA$; see \cite[Chapter~1]{Peller}.

\begin{lemma}\label{lem:bmo-doubling}
If $G_a\in\BMOA$, then $(A_n)$ satisfies \eqref{eq:doubling}.
\end{lemma}
\begin{proof}
Let $x=\one_{[0,2n+2)}$. For $n+1\le m\le2n+1$, at least $n+2$
pairs $(j,k)\in[0,2n+2)^2$ satisfy $j+k=m$. Positivity gives
\[
 (n+2)\sum_{m=n+1}^{2n+1}c_m
 \le\langle\Gamma_{G_a}x,x\rangle
 \le(2n+2)\|\Gamma_{G_a}\|.
\]
Since $G_a\in H^2$, we have $c_m\to0$. Choose $n_*$ so that $c_m\le1/2$
for $m>n_*$. For $n\ge n_*$,
\[
 \log\frac{A_{2n+1}}{A_n}
 =\sum_{m=n+1}^{2n+1}-\log(1-c_m)
 \le2\sum_{m=n+1}^{2n+1}c_m\le4\|\Gamma_{G_a}\|.
\]
The remaining finitely many ratios are finite because $A_n\ge a_0>0$.
\end{proof}

\begin{proposition}[Comparison under doubling]\label{prop:bmo-testing}
If \eqref{eq:doubling} holds with $K\ge2$, then
\begin{equation}\label{eq:bmo-testing}
 \tau(a)\asymp_K1+\mathcal F(G_a).
\end{equation}
\end{proposition}
\begin{proof}
For $J,k\ge1$, positivity and $c_j\le-\log(1-c_j)$ give
\[
 b_k(J)\le\sum_{j=kJ+1}^{(k+1)J}-\log(1-c_j)
 =\log\frac{A_{(k+1)J}}{A_{kJ}}\le\log K,
\]
because $(k+1)J\le2kJ+1$. Hence
\[
 \min\{1,b_k(J)\}^2\le b_k(J)^2
 \le\max\{1,(\log K)^2\}\,\min\{1,b_k(J)\}^2.
\]
Sum over $k$, take the supremum over $J$, and apply
Theorem~\ref{thm:symbol}.
\end{proof}

\begin{proof}[Proof of Corollary~\ref{cor:bmo-intro}]
Under doubling, Proposition~\ref{prop:bmo-testing} and
\eqref{eq:fefferman} give $\tau(a)<\infty$ if and only if
$G_a\in\BMOA$. Apply Theorem~\ref{thm:main}.
Without a prior doubling assumption, $G_a\in\BMOA$ implies both
doubling, by Lemma~\ref{lem:bmo-doubling}, and boundedness, since
$\mathcal F_{\mathrm{tr}}(G_a)\le\mathcal F(G_a)$.
\end{proof}

\subsection{The Bloch condition within the bounded family}

The Bloch space $\mathcal B$ consists of the analytic functions $g$ with
\[
 \|g\|_{\mathcal B,*}
 :=\sup_{z\in\D}(1-|z|^2)|g'(z)|<\infty.
\]
For nonnegative coefficients, its relevance here follows directly
from positivity.

\begin{corollary}\label{cor:bloch-symbol}
Suppose that $N_a$ is bounded on $H^p$ for some $1<p<\infty$. Then
\begin{equation}\label{eq:bloch-symbol}
 G_a\in\mathcal B
 \quad\Longleftrightarrow\quad G_a\in\BMOA
 \quad\Longleftrightarrow\quad (A_n)\text{ is doubling}.
\end{equation}
\end{corollary}

\begin{proof}
First suppose $G_a\in\mathcal B$, and put
$M=\|G_a\|_{\mathcal B,*}$. For an integer $N\ge1$, take $r=e^{-1/N}$.
For $N<j\le2N$, we have $j\ge N$ and $r^{j-1}\ge e^{-2}$, so
\[
 Ne^{-2}\sum_{j=N+1}^{2N}c_j\le G_a'(r)
 \le\frac{M}{1-e^{-2/N}}.
\]
Since $N(1-e^{-2/N})\ge1-e^{-2}$, these block sums are bounded
by $C M$, uniformly in $N$. Taking $N=kJ$ shows that
$b_k(J)\le C M$ for every $k,J\ge1$. Therefore
\[
 \mathcal F(G_a)
 \le\max\{1,C^2M^2\}\,\mathcal F_{\mathrm{tr}}(G_a)<\infty,
\]
where the last inequality follows from boundedness of $N_a$ and
Theorem~\ref{thm:symbol}. By \eqref{eq:fefferman}, $G_a\in\BMOA$.

The implication $G_a\in\BMOA\Rightarrow(A_n)\text{ doubling}$
is Lemma~\ref{lem:bmo-doubling}. Finally, doubling gives
$G_a\in\mathcal B$ directly: for $m\ge0$,
\[
 \sum_{j=2^m+1}^{2^{m+1}}c_j
 \le\log\frac{A_{2^{m+1}}}{A_{2^m}}\le\log K.
\]
Summing the resulting geometric bounds yields
$\sum_{j=1}^{N}jc_j\le C_KN$ for every $N\ge1$.
If $U_N=\sum_{j=1}^{N}jc_j$, summation by parts gives
\[
 G_a'(r)=(1-r)\sum_{N\ge1}U_Nr^{N-1}
 \le\frac{C_K}{1-r}\qquad(0\le r<1).
\]
Positivity implies $|G_a'(z)|\le G_a'(|z|)$, proving the Bloch estimate.
\end{proof}

Thus the truncated criterion describes the full bounded family.
Within this family, the Bloch and BMOA conditions both select exactly
the operators with doubling partial sums.

\section*{Methodology}
For the writing and preparation of this article, the author used the large
language models ChatGPT and Claude, including for proofreading,
bibliographical research and the exploration and checking of arguments. The author remains fully responsible for the content of the paper and
adheres to the principles of the ``Leiden Declaration on Artificial
Intelligence and Mathematics.''

\end{document}